\documentstyle[12pt]{amsart}
\begin{document}
  \title{Bubbling phenomena of  biharmonic maps}
  \title[bubbling phenamena of biharmonic maps]
  {Bubbling phenomena of 
  biharmonic maps}
\author{Nobumitsu Nakauchi}
  \address{Graduate School of Science and Engineering, \newline
  Yamaguchi University, 
  Yamaguchi, 753-8512, Japan}
  \email{nakauchi@@yamaguchi-u.ac.jp}
   \author{Hajime Urakawa}
  \address{Division of Mathematics, Graduate School of Information Sciences, Tohoku University, Aoba 6-3-09, Sendai, 980-8579, Japan}
  \curraddr{Institute for International Education, 
  Tohoku University, Kawauchi 41, Sendai 980-8576, Japan}
  \email{urakawa@@math.is.tohoku.ac.jp}
  \thanks{
  Supported by the Grant-in-Aid for the Scientific Research, (A), No. 19204004, (C) No. 21540207, Japan Society for the Promotion of Science. 
  }
\maketitle
\begin{abstract}
  By using Moser's iteration technique, 
  we show the bubbling phenomena
   for biharmonic maps which cover the ones for harmonic maps. 
  \end{abstract}
\numberwithin{equation}{section}
\theoremstyle{plain}
\newtheorem{df}{Definition}[section]
\newtheorem{th}[df]{Theorem}
\newtheorem{prop}[df]{Proposition}
\newtheorem{lem}[df]{Lemma}
\newtheorem{cor}[df]{Corollary}
\newtheorem{rem}[df]{Remark}
\newtheorem{slem}[df]{Sublemma}
\section{Introduction}
Harmonic maps play a central role in variational problems, which are,  
by definition, critical points of the energy functional 
$E(\varphi)=\frac12\int_M\vert d\varphi\vert^2\,v_g$ 
for smooth maps $\varphi$ of $(M,g)$ into $(N,h)$. 
Extending the notion of harmonic maps, in 1983, J. Eells and L. Lemaire \cite{EL1} proposed the problem 
to consider polyharmonic, i.e., 
$k$-harmonic maps which are critical points of the functional
$$
E_k(\varphi)=\frac12\int_M
\vert(d+\delta)^k\varphi\vert^2\,v_g,\quad (k=1,2,\cdots). 
$$
By definition, we have  $E_1(\varphi)=E(\varphi)$ 
whose Euler-Lagrange equation is that 
the tension field $\tau(\varphi)$ 
vanishes identically 
(cf. \cite{EL1}, \cite{EL2}, \cite{EL3}). 
After G.Y. Jiang \cite{J} studied the first and second variation formulas of $E_2$ given by
$$
E_2(\varphi)=\frac12\int_M\vert\tau(\varphi)\vert^2\,v_g
$$
whose critical points are called (intrinsic) 
{\em biharmonic} maps, there have been extensive studies in this area (for instance, see \cite{CMP}, \cite{LO}, \cite{LO2}, \cite{O1}, \cite{MO1}, \cite{IIU2}, \cite{IIU}. 
\cite{I}, \cite{II}, \cite{S1}, etc.). Harmonic maps are always biharmonic maps by definition. 
\par
The theory of bubbling phenomena of harmonic maps 
was first studied by Sacks and Uhlenbeck \cite{SU} 
and extended to several variational problems 
including Yang-Mills theory (see Freed and Uhlenbeck 
\cite{FU})).  
Then, for the bubbling phenomena of biharmonic maps, we will show
\begin{th} $($cf. {\bf Theorem 3.1}$)$
Let $(M,g)$ and $(N,h)$ be two compact Riemannian manifolds. Assume that $m=\dim M\geq 3$. 
For every positive constant $C>0$, let us consider 
a family of smooth biharmonic maps of $(M,g)$ into $(N,h)$, 
\begin{align}
{\mathcal F}=
\bigg\{
\varphi:&(M,g)\rightarrow (N,h),\,\text{smooth biharmonic}\,\,\vert
\nonumber\\
&\int_M\vert\,d\varphi\,\vert^m\,v_g\leq C
\,\,\text{and}\,\,
\int_M\vert\tau(\varphi)\vert^2\,v_g\leq C
\bigg\}.
\end{align}
Then, any sequence in ${\mathcal F}$ causes a {\em bubbling}: Namely, for any sequence 
$\{\varphi_i\}\subseteq {\mathcal F}$, 
there exist a finite set ${\mathcal S}$ in $M$, say, ${\mathcal S}=\{x_1,\cdots,x_{\ell}\}$, 
and a smooth biharmonic map 
$\varphi_{\infty}:\,(M\backslash{\mathcal S},g)
\rightarrow (N,h)$ 
such that, 
\par
$(1)$ 
a subsequence 
$\varphi_{i_j}$ converges to $\varphi_{\infty}$ in the $C^{\infty}$-topology 
on $M\backslash{\mathcal S}$, as $j\rightarrow\infty$, and 
\par
$(2)$ 
the Radon measures $\vert d\varphi_{i_j}\vert^m \,v_g$ converge to 
a measure 
\begin{equation}
\vert d\varphi_{\infty}\vert^m\,v_g+\sum_{k=1}^{\ell}a_k\,\delta_{x_k}, 
\end{equation}
as $j\rightarrow \infty$.  
Here  $a_k$ is some constant, and 
$\delta_{x_k}$ is the Dirac measure whose support 
is $\{x_k\}$ $(k=1\,\cdots,\ell)$. 
\end{th}
\vskip0.3cm\par
As a corollary, we obtain the bubbling theorem for harmonic maps 
(cf. Theorem 3.2). Theorem 1.1 whose proof will be given in Section 6 
follows from the $C^0$-estimation of tension field $\tau(\varphi)$ (cf. Proposition 4.3) and the $C^1$-estimation 
for biharmonic maps $\varphi$ (cf. Theorem 5.1). 

\vskip0.6cm\par
{\bf Acknowledgement.} 
\quad 
We would like to express our gratitude to Professor H. Naito for his useful suggestions and comments during the preparations of this paper. 
\vskip0.6cm\par
\section{Preliminaries}
In this section, we prepare materials for the first and second variational formulas for the bienergy functional and biharmonic maps. 
Let us recall the definition of a harmonic map $\varphi:\,(M,g)\rightarrow (N,h)$, of a compact Riemannian manifold $(M,g)$ into another Riemannian manifold $(N,h)$, 
which is an extremal 
of the {\em energy functional} defined by 
$$
E(\varphi)=\int_Me(\varphi)\,v_g, 
$$
where $e(\varphi):=\frac12\vert d\varphi\vert^2$ is called the energy density 
of $\varphi$.  
That is, for any variation $\{\varphi_t\}$ of $\varphi$ with 
$\varphi_0=\varphi$, 
\begin{equation}
\frac{d}{dt}\bigg\vert_{t=0}E(\varphi_t)=-\int_Mh(\tau(\varphi),V)v_g=0,
\end{equation}
where $V\in \Gamma(\varphi^{-1}TN)$ is the variation vector field along $\varphi$ which is given by 
$V(x)=\frac{d}{dt}\vert_{t=0}\varphi_t(x)\in T_{\varphi(x)}N$, 
$(x\in M)$, 
and  the {\em tension field} is given by 
$\tau(\varphi)
=\sum_{i=1}^mB(\varphi)(e_i,e_i)\in \Gamma(\varphi^{-1}TN)$, 
where 
$\{e_i\}_{i=1}^m$ is a locally defined frame field on $(M,g)$, 
and $B(\varphi)$ is the second fundamental form of $\varphi$ 
defined by 
\begin{align}
B(\varphi)(X,Y)&=(\widetilde{\nabla}d\varphi)(X,Y)\nonumber\\
&=(\widetilde{\nabla}_Xd\varphi)(Y)\nonumber\\
&=\overline{\nabla}_X(d\varphi(Y))-d\varphi(\nabla_XY)\nonumber\\
&=\nabla^N_{d\varphi(X)}d\varphi(Y)-d\varphi(\nabla_XY),
\end{align}
for all vector fields $X, Y\in {\frak X}(M)$. 
Furthermore, 
$\nabla$, and
$\nabla^N$, 
 are connections on $TM$, $TN$  of $(M,g)$, $(N,h)$, respectively, and 
$\overline{\nabla}$, and $\widetilde{\nabla}$ are the induced ones on $\varphi^{-1}TN$, and $T^{\ast}M\otimes \varphi^{-1}TN$, respectively. By (2.1), $\varphi$ is harmonic if and only if $\tau(\varphi)=0$. 
\par
The second variation formula is given as follows. Assume that 
$\varphi$ is harmonic. 
Then, 
\begin{equation}
\frac{d^2}{dt^2}\bigg\vert_{t=0}E(\varphi_t)
=\int_Mh(J(V),V)v_g, 
\end{equation}
where 
$J$ is an elliptic differential operator, called 
{\em Jacobi operator}  acting on 
$\Gamma(\varphi^{-1}TN)$ given by 
\begin{equation}
J(V)=\overline{\Delta}V-{\mathcal R}(V),
\end{equation}
where 
$\overline{\Delta}V=\overline{\nabla}^{\ast}\overline{\nabla}V
=-\sum_{i=1}^m\{
\overline{\nabla}_{e_i}\overline{\nabla}_{e_i}V-\overline{\nabla}_{\nabla_{e_i}e_i}V
\}$ 
is the {\em rough Laplacian} and 
${\mathcal R}$ is a linear operator on $\Gamma(\varphi^{-1}TN)$
given by 
${\mathcal R}(V)=
\sum_{i=1}^mR^N(V,d\varphi(e_i))d\varphi(e_i)$,
and $R^N$ is the curvature tensor of $(N,h)$ given by 
$R^N(U,V)=\nabla^N{}_U\nabla^N{}_V-\nabla^N{}_V\nabla^N{}_U-\nabla^N{}_{[U,V]}$ for $U,\,V\in {\frak X}(N)$.   
\par
J. Eells and L. Lemaire \cite{EL1} proposed polyharmonic ($k$-harmonic) maps and 
Jiang \cite{J} studied the first and second variation formulas of biharmonic maps. Let us consider the {\em bienergy functional} 
defined by 
\begin{equation}
E_2(\varphi)=\frac12\int_M\vert\tau(\varphi)\vert ^2v_g, 
\end{equation}
where 
$\vert V\vert^2=h(V,V)$, $V\in \Gamma(\varphi^{-1}TN)$.  
\par
Then, the first variation formula of the bienergy functional 
is given as follows. 
\begin{th}
\quad $($the first variation formula$)$ 
\begin{equation}
\frac{d}{dt}\bigg\vert_{t=0}E_2(\varphi_t)
=-\int_Mh(\tau_2(\varphi),V)v_g.
\end{equation}
Here, 
\begin{equation}
\tau_2(\varphi)
:=J(\tau(\varphi))=\overline{\Delta}\tau(\varphi)-{\mathcal R}(\tau(\varphi)),
\end{equation}
which is called the {\em bitension field} of $\varphi$, and 
$J$ is given in $(2.4)$.  
\end{th}
\begin{df}
A smooth map $\varphi$ of $M$ into $N$ is said to be 
{\em biharmonic} if 
$\tau_2(\varphi)=0$. 
\end{df} 
\vskip0.6cm\par
\section{Bubbling Theorem of Biharmonic Maps}
We have the following bubbling theorem for biharmonic maps. 
\begin{th} $($Bubbling for Biharmonic Maps$)$
Let $(M,g)$ and $(N,h)$ be two compact Riemannian manifolds without boundaries. Assume that $\dim M=m\geq 3$. 
For every positive constant $C>0$, consider 
a family of smooth 
biharmonic maps of $(M,g)$ into $(N,h)$, 
\begin{align}
{\mathcal F}=
\bigg\{
\varphi:&(M,g)\rightarrow (N,h),\,\text{smooth biharmonic}\,\,\vert
\nonumber\\
&\int_M\vert\,d\varphi\,\vert^m\,v_g\leq C
\,\,\text{and}\,\,
\int_M\vert\tau(\varphi)\vert^2\,v_g\leq C
\bigg\}.
\end{align}
Then, any sequence in ${\mathcal F}$ causes a {\bf bubbling}: Namely, for any sequence 
$\{\varphi_i\}\subseteq {\mathcal F}$, 
there exist a finite set ${\mathcal S}$ in $M$, say, ${\mathcal S}=\{x_1,\cdots,x_{\ell}\}$, 
and a smooth biharmonic map 
$\varphi_{\infty}:\,(M\backslash{\mathcal S},g)
\rightarrow (N,h)$ 
such that, 
\par
$(1)$ 
a subsequence 
$\varphi_{i_j}$ converges to $\varphi_{\infty}$ in the $C^{\infty}$-topology 
on $M\backslash{\mathcal S}$, as $j\rightarrow\infty$, and 
\par
$(2)$ 
the Radon measure $\vert d\varphi_{i_j}\vert^m \,v_g$ converges to 
a measure 
\begin{equation}
\vert d\varphi_{\infty}\vert^m\,v_g+\sum_{k=1}^{\ell}a_k\,\delta_{x_k}, 
\end{equation}
as $j\rightarrow \infty$.  
Here  $a_k$ is some constant, and 
$\delta_{x_k}$ is the Dirac measure centered at $\{x_k\}$ $(k=1\,\cdots,\ell)$. 
\end{th}
\vskip0.6cm\par
As a corollary, we have immediately 
\begin{th} $($Bubbling for Harmonic Maps$)$
Let $(M,g)$ and $(N,h)$ be two compact Riemannian manifolds without boundaries. 
Assume that $\dim M=m\geq 3$. 
For every positive constant $C>0$, let us consider 
a family of smooth harmonic maps of $(M,g)$ into $(N,h)$, 
\begin{equation}
{\mathcal F}^h=\left\{
\varphi:(M,g)\rightarrow (N,h),\,\text{smooth harmonic}\,\,\vert
\int_M\vert\,d\varphi\,\vert^m\,v_g\leq C
\right\}.
\end{equation}
Then, any sequence in ${\mathcal F}^h$ causes a {\bf bubbling}: Namely, for any sequence 
$\{\varphi_i\}\in {\mathcal F}^h$, 
there exist a finite set ${\mathcal S}$ in $M$, say, ${\mathcal S}=\{x_1,\cdots,x_{\ell}\}$, 
and a smooth harmonic map 
$\varphi_{\infty}:\,(M\backslash{\mathcal S},g)
\rightarrow (N,h)$ 
such that, 
\par
$(1)$ 
a subsequence 
$\varphi_{i_j}$ converges to $\varphi_{\infty}$ in the $C^{\infty}$-topology 
on $M\backslash{\mathcal S}$, as $j\rightarrow\infty$, and 
\par
$(2)$ 
the Radon measures $\vert d\varphi_{i_j}\vert^m \,v_g$ converge to 
a measure 
\begin{equation}
\vert d\varphi_{\infty}\vert^m\,v_g+\sum_{k=1}^{\ell}a_k\,\delta_{x_k}, 
\end{equation}
as $j\rightarrow \infty$.  
\end{th}
\begin{pf} 
For any sequence in $\{\varphi_i\}\in {\mathcal F}^h$, 
the limit $\varphi_{\infty}$ in Theorem 3.1 
is a smooth biharmonic map of 
$(M\backslash {\mathcal S},g)$ into $(N,h)$. 
Due to $(1)$ of Theorem 3.1, 
$\varphi_{i_j}$ converges to 
$\varphi_{\infty}$ in the $C^{\infty}$-topology on 
$M\backslash{\mathcal S}$, 
so that $\tau(\varphi_{i_j})$ converges to 
$\tau(\varphi_{\infty})$ pointwise 
on $M\backslash{\mathcal S}$. 
Since $\tau(\varphi_{i_j})\equiv0$, we have 
$\tau(\varphi_{\infty})\equiv0$ on 
$M\backslash{\mathcal S}$, i.e., 
$\varphi_{\infty}$ is harmonic on $M\backslash{\mathcal S}$. 
And $(1)$ and $(2)$ hold also due to Theorem 3.1. 
\end{pf}
\vskip0.6cm\par
\section{Basic Inequalities}
To prove Theorem 3.1, it is necessary to prepare
 the following two basic inequalities. 
\begin{lem}
 Assume that the sectional curvature of $(N,h)$ is bounded 
 above by a constant $C$, 
 and $\varphi:\,(M\backslash {\mathcal S},g)\rightarrow (N,h)$ is a biharmonic mapping for some closed set ${\mathcal S}$ of $M$. 
 Then, 
 we have 
 \begin{equation}
 \frac12\,\Delta\,\vert \tau(\varphi)\vert^2+C\,\vert d\varphi\,\vert^2
 \,\vert \tau(\varphi)\vert^2\geq 
 \vert\overline{\nabla}\tau(\varphi)\vert^2.
 \end{equation}
\end{lem}
\begin{pf}
\quad 
Let us take a local orthonormal frame field 
$\{e_i\}_{i=1}^m$ 
on $M\backslash{\mathcal S}$, and $\varphi:\,(M\backslash{\mathcal S},g)\rightarrow (N,h)$, a biharmonic map. 
Then, for $V:=\tau(\varphi)\in \Gamma(\varphi^{-1}TN)$, we have 
\begin{align}
\frac12\,\Delta\,\vert V\vert^2
&=\frac12
\sum_{i=1}^m
\left\{
e_i{}^2\,\vert V\vert^2-\nabla_{e_i}e_i\,\vert V\vert^2
\right\}\nonumber\\
&=
\sum_{i=1}^m
\left\{
e_i\,h(\overline{\nabla}_{e_i}V,V)
-h(\overline{\nabla}_{\nabla_{e_i}e_i}V,V)
\right\}\nonumber\\
&=
\sum_{i=1}^m
\left\{
h(\overline{\nabla}_{e_i}\overline{\nabla}_{e_i}\,V,V)
-h(\overline{\nabla}_{\nabla_{e_i}e_i}\,V,V)
\right\}
\nonumber\\
&\qquad
+\sum_{i=1}^mh(\overline{\nabla}_{e_i}V,
\overline{\nabla}_{e_i}V)
\nonumber\\
&=h(-\overline{\Delta}\,V,V)
+\vert\overline{\nabla}V\vert^2\nonumber\\
&=h(-{\mathcal R}(V),V)+\vert\overline{\nabla}V\vert^2, 
\end{align}
because for the second last equality, we used $\overline{\Delta}V-{\mathcal R}(V)=J(V)=0$ for $V=\tau(\varphi)$, 
due to the biharmonicity of $\varphi:\,(M\backslash {\mathcal S},g)\rightarrow (N,h)$. 
Since 
$$h({\mathcal R}(V),V)=
\sum_{i=1}^mh(R^N(V,d\varphi(e_i))
d\varphi(e_i),V),$$ 
the right hand side of (4.2) is bigger than or equal to 
$$-C\sum_{i=1}^m\vert d\varphi(e_i)\vert^2\,\vert V\vert^2
+\vert\,\overline{\nabla}V\,\vert^2
=-C\,\vert d\varphi\vert^2\,\vert V\vert^2
+\vert\,\overline{\nabla}V\,\vert^2.$$
We have (4.1). 
\end{pf}
\vskip0.6cm\par
\begin{lem}
Under the same assumption as Lemma 4.1, we have 
\begin{equation}
\vert\tau(\varphi)\vert\,\Delta\,\vert\tau(\varphi)\vert
+C\,\vert d\varphi\,\vert^2\,\vert\tau(\varphi)\vert^2\geq 0
\end{equation}
for all $\varphi\in C^{\infty}(M,N)$. 
\end{lem}
\begin{pf}
\quad
By the following equality 
\begin{align}
\Delta\,\vert\tau(\varphi)\vert^2
=2\vert\tau(\varphi)\vert\,\Delta\vert\tau(\varphi)\vert
+\vert\nabla\vert\tau(\varphi)\vert\,\vert^2,
\end{align}
and Lemma 4.1, 
we have 
\begin{align}
\vert\tau(\varphi)\vert\,\Delta\,\vert\tau(\varphi)\vert
&+\vert\,\nabla\,\vert\tau(\varphi)\vert\,\vert^2
+C\,\vert d\varphi\,\vert^2\,\vert\tau(\varphi)\vert^2
\nonumber\\
&\geq\frac12\Delta\,\vert\tau(\varphi)\vert
+C\,\vert d\varphi\,\vert^2\,\vert\tau(\varphi)\vert^2\nonumber\\
&\geq  \vert\overline{\nabla}\tau(\varphi)\vert^2.
\end{align}
So that we have 
\begin{align}
\vert\tau(\varphi)\vert\,\Delta\,\vert\tau(\varphi)\vert
&+C\,\vert d\varphi\,\vert^2\,\vert\tau(\varphi)\vert^2\nonumber\\
&\geq 
 \vert\overline{\nabla}\tau(\varphi)\vert^2
 -\vert\,\nabla\,\vert\tau(\varphi)\vert\,\vert^2\nonumber\\
 &\geq 0. 
\end{align}
Here, to see the last inequality of (4.6), it suffices to 
notice that for all $V\in \Gamma(\varphi^{-1}TN)$, 
\begin{equation}
\vert\overline{\nabla}V\vert\geq \vert\,\,\nabla\,\vert V\vert\,\,\vert 
\end{equation}
which follows from that 
\begin{align}
\vert V\vert\,\,\vert\,\nabla\,\vert V\vert\,\,\vert
&=\frac12\,\vert\,\nabla\,\vert V\,\vert^2\,\,\vert\nonumber\\
&=\frac12\,\vert \nabla\,h(V,V)\,\vert\nonumber\\
&=\vert\,h(\overline{\nabla}V,V)\,\vert\nonumber\\
&\leq \vert\,\overline{\nabla}V\,\vert\,\,\vert V\vert.
\end{align}
This proves Lemma 4.2. 
\end{pf}
\vskip0.6cm\par
The following proposition is to give the 
$C^0$-estimation for the tension field of a biharmonic map which is necessary together with Theorem 5.1 in order to give 
the $C^1$-estimation of $\varphi$. 
\begin{prop}
Assume that 
the sectional curvature of $(N,h)$ is bounded above by 
a positive constant $C>0$. 
Then, there exists a positive number 
$\varepsilon_0$ depending only on 
the Sobolev constant of $(M,g)$ and $C$ such that 
for every smooth biharmonic map 
$\varphi$ of $(M,g)$ into $(N,h)$, if 
\begin{equation}
\int_{B_{2r}(x_0)}\vert d\varphi\vert^m\,v_g\leq \varepsilon_0, 
\end{equation}
then 
\begin{equation}
\sup_{B_{r}(x_0)}\vert\tau(\varphi)\vert^2\leq
\frac{C'}{r^{m/2}}\int_{B_{2r}(x_0)}\vert\tau(\varphi)\vert^2\,v_g.
\end{equation} 
for some positive constant $C'>0$ depending only on $C$ and $m=\dim M$. 
\end{prop}
\vskip0.6cm\par
\begin{pf}
({\it The first step}) \quad 
For every 
$0<\rho_1<\rho_2\leq 2r<\infty$, 
we first take a cutoff  $C^{\infty}$ function $\eta$ on $M$ 
(for instance, see \cite{K}) satisfying that 
\begin{equation}
\left\{
\begin{aligned}
0\leq &\eta(x)\leq 1\quad (x\in M),\\
\eta(x)&=1\qquad\quad (x\in B_{\rho_1}(x_0)),\\
\eta(x)&=0\qquad\quad (x\not\in B_{\rho_2}(x_0)),\\
\vert\nabla\eta\vert&\leq\frac{2}{\rho_2-\rho_1}
\quad (x\in M).
\end{aligned}
\right.
\end{equation}
\par
Multiply $\vert\tau(\varphi)\vert^{p-2}\,\eta^2$ to both hand sides of $(4.3)$ 
and integrate it over $M$. 
Then, we have 
\begin{equation}
0\leq 
\int_M\vert\tau(\varphi)\vert^{p-1}\,\eta^2\,\Delta(\vert\tau(\varphi)\vert)\,v_g
+C\int_M\vert d\varphi\vert^2\,\vert\tau(\varphi)\vert^p\,\eta^2\,v_g.
\end{equation}
\par
In order to estimate 
the second term of the right hand side of (4.8), 
we need the following lemma whose proof is omitted since it can be proved immediately 
by H\"older's inequality and 
Sobolev-Poincar\'{e} inequality (see (4.20) and (4.21) below).
\vskip0.6cm\par
\begin{lem}
We have 
\begin{align}
\int_M&\vert d\varphi\vert^2\,\vert\tau(\varphi)\vert^p\,\eta^2\,v_g\nonumber\\
&\leq
C''\,\left\{
\int_{B_{\rho_2}(x_0)}\vert d\varphi\vert^m\,v_g
\right\}^{2/m}
\,\int_M
\vert\nabla(\,\vert\tau(\varphi)\vert^{p/2}\,\eta)\,\vert^2\,v_g,
\end{align}
where $C''>0$ is a positive constant independent on 
$\varphi$. 
\end{lem}
\vskip0.6cm\par
We continue 
the proof of Proposition 4.3. 
We first obtain by (4.12), 
\begin{align}
0&\leq 
\int_M\vert\tau(\varphi)\vert^{p-1}\,\eta^2\,\Delta(\vert\tau(\varphi)\vert)\,v_g
+C\int_M\vert d\varphi\vert^2\,\vert\tau(\varphi)\vert^p\,\eta^2\,v_g
\nonumber\\
&
=
-\frac{4(p-1)}{p^2}\int_M\vert \nabla(\vert\tau(\varphi)\vert^{p/2})\vert^2\,\eta^2\,v_g\nonumber\\
&\quad-\frac{4}{p}\int_M
g(\eta\,\nabla(\vert\tau(\varphi)\vert^{p/2}),
\vert\tau(\varphi)\vert^{p/2}\,\nabla\eta)
\,v_g\nonumber\\
&\quad
+C\int_M\vert d\varphi\vert^2\,\vert\tau(\varphi)\vert^p\,\eta^2\,v_g.
\end{align}
Then, we obtain the following: 
\begin{align}
\int_M\eta^2\,\vert\,\nabla(\vert\tau(\varphi)\vert^{p/2})\vert^2\,v_g
&\leq \frac{p^2}{(p-1)^2}\,
\int_M
\vert\tau(\varphi)\vert^p\,\vert\nabla\eta\vert^2\,v_g\nonumber\\
&\quad+C\int_M\vert d\varphi\vert^2\,\vert\tau(\varphi)\vert^p\,\eta^2\,v_g.
\end{align} 
And then, we have, 
\begin{align}
&\int_M\vert\,\nabla(\vert\tau(\varphi)\vert^{p/2}\,\eta)\,\vert^2\,v_g
\nonumber\\
&\leq
4\,\frac{p^2}{(p-1)^2}\int_M\vert\tau(\varphi)\vert^p\,\vert\nabla\eta\vert^2\,v_g
+C\int_M\vert d\varphi\vert^2\,\vert\tau(\varphi)\vert^p\,\eta^2\,v_g\nonumber\\
&\leq
4\,\frac{p^2}{(p-1)^2}\int_M\vert\tau(\varphi)\vert^p\,\vert\nabla\eta\vert^2\,v_g\nonumber\\
&\quad +CC''\,\left\{
\int_{B_{\rho_2}(x_0)}\vert d\varphi\vert^m\,v_g
\right\}^{2/m}
\,\int_M
\vert\nabla(\,\vert\tau(\varphi)\vert^{p/2}\,\eta)\,\vert^2\,v_g.
\end{align}
In the last inequality, we used (4.13) 
in Lemma 4.4. 
\par
Assume that 
\begin{equation}
\int_{B_{\rho_2}(x_0)}\vert d\varphi\vert^m\,v_g\leq \varepsilon_0. 
\end{equation}
Then, due to (4.16), we have 
\begin{align}
\int_M\vert\,\nabla(\vert\tau(\varphi)\vert^{p/2}\,\eta)\,\vert^2\,v_g
&\leq
4\,\frac{p^2}{(p-1)^2}\int_M\vert\tau(\varphi)\vert^p\,\vert\nabla\eta\vert^2\,v_g\nonumber\\
&\quad+CC''\varepsilon_0{}^{2/m}\,\int_M
\vert\nabla(\,\vert\tau(\varphi)\vert^{p/2}\,\eta)\,\vert^2\,v_g.
\end{align}
If we take $\varepsilon_0>0$ small enough such that 
$1-CC''\,\varepsilon_0{}^{2/m}>\frac12$, i.e., 
$\frac{1}{(2CC'')^{2/m}}>\varepsilon_0$, then, 
we have 
\begin{align}
\int_M\vert\,\nabla(\vert\tau(\varphi)\vert^{p/2}\,\eta)\,\vert^2\,v_g
&\leq
\frac{8p^2}{(p-1)^2}\int_M\vert\tau(\varphi)\vert^p\,\vert\nabla\eta\vert^2\,v_g\nonumber\\
&\leq \frac{8p^2}{(p-1)^2}\frac{4}{(\rho_2-\rho_1)^2}
\int_{B_{\rho_2}(x_0)}\vert\tau(\varphi)\vert^p\,v_g.
\end{align}
\par
For the left hand side of (4.19), 
let us recall the Sobolev embedding theorem (cf. \cite{A}, p. 55; 
 \cite{FU}, p. 95): 
\begin{equation}
H^2_1(M)\subset L^{\gamma}(M),
\end{equation}
where $\gamma:=\frac{2m}{m-2}$, i.e., 
there exists a pocitive constant $C>0$ such that 
\begin{equation}
\left(
\int_M
\vert\,f\,\vert^{\gamma}\,v_g
\right)^{1/\gamma}
\leq
C\left(
\int_M\vert\,\nabla f\,\vert^2\,v_g
\right)^{1/2}
\qquad (\forall \,f\in H^2_1(M)).
\end{equation}
\vskip0.3cm\par
Therefore, we have 
\begin{align}
\int_M\vert\,\nabla(\vert\tau(\varphi)\vert^{p/2}\,\eta)\,\vert^2\,v_g
&\geq 
\frac{1}{C}
\left(
\int_M
\left\{
\,\vert\tau(\varphi)\vert^{p/2}\,\eta
\right\}^{\gamma}
\,v_g
\right)^{2/\gamma}
\nonumber\\
&\geq
\frac{1}{C}
\left(
\int_{B_{\rho_1}(x_i)}
\left(\vert\tau(\varphi)\vert^{p/2}\right)^{\gamma}\,v_g
\right)^{2/\gamma}.
\end{align}
Thus, together with (4.19) and (4.22), we have 
\begin{lem} 
Assume that $(M,g)$ is a compact Riemannian manifold, 
and $\varphi:\,(M\backslash\,x_0,g)\rightarrow (N,h)$ 
is a biharmonic mapping. 
Then, 
for each $0<\rho_1<\rho_2\leq 2r<\infty$, and $2\leq p<\infty$, 
it holds that for each $i=1,\cdots,\ell$, 
\begin{align}
\left(
\int_{B_{\rho_1}(x_0)}
\left(\vert\tau(\varphi)\vert^{p/2}\right)^{\gamma}\,v_g
\right)^{1/\gamma}
&\leq
\frac{p}{p-1}\,\,
\frac{C'}{\rho_2-\rho_1}\nonumber\\
&\qquad
\times
\left(\int_{B_{\rho_2}(x_0)}
\left(\vert\tau(\varphi)\vert^{p/2}\right)^2\,v_g
\right)^{1/2},
\end{align}
where $C'=4\,\sqrt{C}$, and 
$C>0$ is the Sobolev constant in $(4.21)$ 
and $\gamma:=\frac{2m}{m-2}$, $m=\dim M\geq 3$. 
\end{lem}
\vskip0.6cm\par
Now, the Moser's iteration works well, and then 
we obtain Proposition 4.3. 
Indeed, the procedure goes as follows. 
\vskip0.6cm\par
({\it The second step}) \quad 
Here, let us define 
\begin{equation}
\left\{
\begin{aligned}
\overline{\gamma}&:=\frac{m}{m-2}=\frac12\gamma,\\
p_k&:=2\,\overline{\gamma}^{k-1}\,\,\rightarrow\infty\quad (k\rightarrow\infty),\\
r_k&:=\left(1+\frac{1}{2^{k-1}}\right)\,r
\,\,\rightarrow r\quad(k\rightarrow\infty), 
\end{aligned}
\right.
\end{equation}
and in (4.23), let us put 
\begin{equation}
\left\{
\begin{aligned}
p&:=p_k,\\
\rho_1&:=r_{k+1},\\
\rho_2&:=r_k.
\end{aligned}
\right.\nonumber
\end{equation}
Then, we have 
\begin{equation}
\left\{
\begin{aligned}
\frac{p_k\,\gamma}{2}&=p_k\,\overline{\gamma}=2\,\overline{\gamma}^k=p_{k+1},\\
\rho_2-\rho_1&=r_k-r_{k+1}=\left(
\frac{1}{2^{k-1}}-\frac{1}{2^k}
\right)\,r=\frac{1}{2^k}\,r,
\end{aligned}
\right.
\end{equation}
so that (4.23) can be rewritten as follows. 
\begin{align}
\left(\int_{B_{r_{k+1}}(x_i)}\vert\tau(\varphi)\vert^{p_{k+1}}\,v_g
\right)^{1/\gamma}
&\leq
\frac{2\,\overline{\gamma}^{k-1}}{2\,\overline{\gamma}^{k-1}-1}
\,
\,\frac{2^k}{r}\nonumber\\
&\qquad\times
\left(
\int_{B_{r_k}(x_i)}
\vert\tau(\varphi)\vert^{p_k}\,v_g
\right)^{1/2}.
\end{align}
By taking $\frac{1}{\overline{\gamma}^{k-1}}$ power of (4.26), 
we have 
\begin{align}
\Vert
\tau(\varphi)
\Vert_{L^{p_{k+1}}(B_{r_{k+1}}(x_i))}
&\leq 
\left(
\frac{2\,\overline{\gamma}^{k-1}}{2\,\overline{\gamma}^{k-1}-1}
\right)^{2/p_k}
\,\frac{2^{(k/\overline{\gamma}^{k-1})}}{r^{(1/\overline{\gamma}^{k-1})}}
\nonumber\\
&\qquad
\times\,\Vert\tau(\varphi)\Vert_{L^{p_k}(B_{r_k}(x_i))}
\end{align}
since we calculate the power of the left hand side of (4.26) as 
$$
\frac{1}{\gamma}\,\frac{1}{\overline{\gamma}^{k-1}}
=\frac{1}{2\,\overline{\gamma}\,\,\overline{\gamma}^{k-1}}
=\frac{1}{2\,\overline{\gamma}^k}=\frac{1}{p_{k+1}}.
$$
\vskip0.6cm\par
({\it The third step}) \quad 
Now iterate (4.27), then we have 
 \begin{align}
 \Vert\tau(\varphi)\Vert_{L^{p_{k+1}}(B_{r_{k+1}}(x_i))}
 &\leq
 \prod_{k=1}^{\infty}
 \left(
 \frac{2\,\overline{\gamma}^{k-1}}{2\,\overline{\gamma}^{k-1}-1}
 \right)^{2/p_k}
 \,\,\frac{2^{(k/\overline{\gamma}^{k-1})}}{r^{(1/\overline{\gamma}^{k-1})}}
 \nonumber\\
 &\qquad\times
 \Vert\tau(\varphi)\Vert_{L^2(B_{2r}(x_i))}
 \end{align}
 since 
 $p_1=2$ and $r_1=2r$. Here, we notice that 
 \begin{align}
 \prod^{\infty}_{k=1}\frac{1}{r^{(1/\overline{\gamma}^{k-1})}}
 =\frac{1}{r^{(\sum_{k=1}^{\infty} 1/\overline{\gamma}^{k-1})}}
 =\frac{1}{r^{m/2}}
 \end{align}
 since 
 $$
 \sum_{k=1}^{\infty}\frac{1}{\overline{\gamma}^{k-1}}
 =\frac{1}{1-\frac{1}{\overline{\gamma}}}=
 \frac{1}{1-\frac{m-2}{m}}=\frac{m}{2}.
 $$
 Notice also that 
 \begin{equation}
 \prod_{k=1}^{\infty}
 \frac{1}{(2\,\overline{\gamma}^{k-1}-1)^{2/p_k}}\leq1
 \end{equation}
 since 
 $2\,\overline{\gamma}^{k-1}-1> 2-1=1$ 
 when $\overline{\gamma}=m/(m-2)> 1$ $(m\geq 3)$. 
And also notice that 
\begin{align}
\prod_{k=1}^{\infty}2^{(k/\overline{\gamma}^{k-1})}&
=2^{\sum_{k=1}^{\infty}\frac{k}{\overline{\gamma}^{k-1}}}<\infty,\\
\prod_{k=1}^{\infty}
\left(2\,\overline{\gamma}
\right)^{2(k-1)/p_k}
=\gamma&^{2\,\sum_{k=1}^{\infty}\frac{k-1}{p_k}}
=\gamma^{\sum_{k=1}^{\infty}\frac{k-1}{\overline{\gamma}^{k-1}}}<\infty.
\end{align}
Therefore, (4.28) can be written as 
\begin{equation}
\Vert\tau(\varphi)\Vert_{L^{p_{k+1}}(B_{r_{k+1}}(x_i))}
 \leq
 C''\,\frac{1}{r^{m/2}}\,\Vert\tau(\varphi)\Vert_{L^2(B_{2r}(x_i))}
\end{equation}
for some positive constant $C''$ depending only on $m=\dim M$. 
\vskip0.6cm\par
({\it The fourth step}) \quad 
Now, let $k$ tend to infinity. Then, by (4.24), the norm 
$\Vert\tau(\varphi)\Vert_{L^{p_{k+1}}(B_{r_{k+1}}(x_i))}$ 
tends to 
$$
\Vert\tau(\varphi)\Vert_{L^{\infty}(B_r(x_i))}
=\sup_{B_r(x_i)}\vert\tau(\varphi)\vert.
$$
Thus, we obtain 
\begin{equation}
\sup_{B_r(x_i)}\vert\tau(\varphi)\vert\leq
\frac{C''}{r^{m/2}}\,\Vert\tau(\varphi)\Vert_{L^2(B_{2r}(x_i))},
\end{equation}
which is the desired inequality (4.10). 
We have Proposition 4.3. 
\end{pf}
\vskip0.6cm\par
\section{$C^1$-estimate for biharmonic maps}
\subsection{Statement of the results}
In this section, we will show the following theorem which gives the $C^1$-estimation of biharmonic maps together with 
Proposition 4.3: 
\begin{th}
Assume that $(M,g)$ is a compact 
Riemannian manifold without boundary 
and 
the sectional curvature of $(N,h)$ 
is bounded from above. 
Then, there exist three positive constants
$\varepsilon_1>0$, $\varepsilon_2>0$  and $C^{\ast}>0$ 
which depend only on the geometry 
of $(M,g)$ such that,  
if $\varphi:\,(M,g)\rightarrow (N,h)$ 
is biharmonic with finite bienergy 
$E_2(\varphi)<\infty$, 
and 
satisfies that 
\begin{equation}
\int_{B_r(x_0)}\vert d\varphi\vert^m\,v_g<\varepsilon_1
\end{equation}
and 
\begin{equation}
\int_{B_r(x_0)}\vert \tau(\varphi)\vert^m\,v_g
<\varepsilon_2,
\end{equation}
then, 
 \begin{align}
\sup_{B_{r/2}(x_0)}\vert d \varphi\vert 
+\sup_{B_{r/2}(x_0)}\vert\tau(\varphi)\vert
\leq\frac{C^{\ast}}{r}
\bigg[
\varepsilon_1{}^{1/m}
+
\varepsilon_2{}^{1/m}
+1
\bigg]
\end{align}
where $m=\dim M$. 
\end{th}
\vskip0.6cm\par
\par
Combining Theorem 5.1 and Proposition 4.3, we will have the $C^1$-estimate for a biharmonic map with finite bienergy, which will be used 
in the proof of Theorem 3.1.  
\subsection{Sublemmas on Bochner-type estimate}
To prove it, we need three Lemmas concerning on 
the $L^{1,q}$-estimates of $\vert d\varphi\vert$ of a bi-harmonic map 
$\varphi:\,(M,g)\rightarrow (N,h)$.   
To obtain them, we first prepare two sublemmas 
on Bochner-type formulas for a bi-harmonic mapping. 
 \begin{slem} 
 Assume that $(M,g)$ is an $m$-dimensional Riemannian manifold, and the sectional curvature of 
 $(N,h)$ is bounded from above, then,  
 for every smooth map $\varphi$ of $M$ into $N$, 
 it holds that 
 \begin{align}
 \frac12 \Delta\,& \vert d\varphi\vert^2
 = \vert\widetilde{\nabla} d\varphi\vert^2
 +\sum_{i=1}^mh(\overline{\nabla}_{e_i}\tau(\varphi),d\varphi(e_i))\nonumber\\
 &-\sum_{i,j=1}^mh({}^N\!R(d\varphi(e_i),d\varphi(e_j))d\varphi(e_j),d\varphi(e_i))\nonumber\\
 &-\sum_{i,j=1}^m\{
 h((\widetilde{\nabla}_{e_j}d\varphi)(\nabla_{e_i}e_j),d\varphi(e_i))
 +h((\widetilde{\nabla}_{\nabla_{e_i}e_j}d\varphi)(e_j),d\varphi(e_i))
 \}, 
 \end{align} 
 where the Laplacian $\Delta$ on 
 $C^{\infty}(M)$ is defined by 
 $
 \Delta f=\sum_{i=1}^m\{e_i(e_i f)-\nabla_{e_i}e_if\}, \quad (f\in C^{\infty}(M))
 $. 
 \end{slem}
 \begin{pf} 
 For a completeness, we give a proof. 
 We may assume $\nabla e_i=0$ at every fixed point $x\in M$. 
 Then, 
 for the left hand side of (5.4),  we first notice that  
 \begin{align}
 \frac12\Delta\vert d\varphi\vert^2
 &= \sum_{j=1}^m \langle
  \widetilde{\nabla}_{e_j}(\widetilde{\nabla}_{e_j}d\varphi)
  -\widetilde{\nabla}_{{\nabla}_{e_j}e_j} d\varphi,d\varphi\rangle
  +\sum_{j=1}^m
  \langle \widetilde{\nabla}_{e_j}d\varphi,
  \widetilde{\nabla}_{e_j}d\varphi\rangle
 \nonumber\\
 &=\sum_{i,j=1}^m\langle((\widetilde{\nabla}_{e_j}\widetilde{\nabla}_{e_j}d\varphi)(e_i),d\varphi(e_i)\rangle+\vert \widetilde{\nabla} d\varphi\vert^2.
 \end{align}
 Indeed, 
 the second equality of (5.5) follows from that, 
 for each vector fields $X, Y, Z$ 
on $M$, 
 \begin{align}
 (\widetilde{\nabla}_X\widetilde{\nabla}_Y
 d\varphi)(Z)
 &=
 \overline{\nabla}_X((\widetilde{\nabla}_Yd\varphi)(Z))-
 (\widetilde{\nabla}_{\nabla_XY}d\varphi)(Z)
 -(\widetilde{\nabla}_Yd\varphi)(\nabla_XZ)\nonumber\\
 &=\widetilde{\nabla}_X(\widetilde{\nabla}_Y
 d\varphi)(Z)-
 (\widetilde{\nabla}_{\nabla_XY}d\varphi)(Z). 
  \end{align}
  \par
  Then, for the first term of (5.5), 
  we have 
  \begin{align}
  (\widetilde{\nabla}_{e_j}\widetilde{\nabla}_{e_j}d\varphi)(e_i)&
  =(\widetilde{\nabla}_{e_j}\widetilde{\nabla}_{e_i}
  d\varphi)(e_j)\nonumber\\
  &
  =(\widetilde{\nabla}_{e_i}\widetilde{\nabla}_{e_j}
  d\varphi)(e_j)
  +(R^{\widetilde{\nabla}}(e_j,e_i)d\varphi)(e_j)
  \nonumber\\
  &
  =(\widetilde{\nabla}_{e_i}\widetilde{\nabla}_{e_j}
  d\varphi)(e_j)
  +{}^N\!R(d\varphi(e_j),d\varphi(e_i))d\varphi(e_j)
  \nonumber\\
  &=
  (\widetilde{\nabla}_{e_i}\widetilde{\nabla}_{e_j}
  d\varphi)(e_j)
  -{}^N\!R(d\varphi(e_i),d\varphi(e_j))d\varphi(e_j).
  \end{align} 
  Here, 
  the first equality of (5.7) 
  follows from that:  
  the left hand side of (5.7) is equal to 
  \begin{align}
  \overline{\nabla}_{e_j}&((\widetilde{\nabla}_{e_j}d\varphi)(e_i))
  -(\widetilde{\nabla}_{e_j}d\varphi)(\nabla_{e_j}e_i)-(\widetilde{\nabla}_{\nabla_{e_j}e_j}d\varphi)(e_i)\nonumber\\
  &=\overline{\nabla}_{e_j}((\widetilde{\nabla}_{e_j}d\varphi)(e_i))
\nonumber\\
&=\overline{\nabla}_{e_j}
(\overline{\nabla}_{e_j}
(d\varphi(e_i)))
-\overline{\nabla}_{e_j}(d\varphi(\nabla_{e_j}e_i)), 
  \end{align}
at $x\in M$, and by the same manner, the right hand side of the first line of (5.7) is equal to 
 \begin{equation}
  \overline{\nabla}_{e_j}
(\overline{\nabla}_{e_i}
(d\varphi(e_j)))
-\overline{\nabla}_{e_j}(d\varphi(\nabla_{e_i}e_j))
\end{equation}
at $x\in M$. Then, at $x\in M$, 
$(\widetilde{\nabla}_{e_j}\widetilde{\nabla}_{e_j}d\varphi)(e_i)-
(\widetilde{\nabla}_{e_j}\widetilde{\nabla}_{e_i}d\varphi)(e_j)$ is equal to
\begin{align}
\overline{\nabla}_{e_j}
&\{\overline{\nabla}_{e_j}(d\varphi(e_i))
-\overline{\nabla}_{e_i}(d\varphi(e_j))\}
-\overline{\nabla}_{e_j}(d\varphi(\nabla_{e_j}e_i))
+
\overline{\nabla}_{e_j}(d\varphi(\nabla_{e_j}e_i))\nonumber\\
&=\overline{\nabla}_{e_j}\{d\varphi([e_j,e_i])\}
-\overline{\nabla}_{e_j}\{
d\varphi([e_j,e_i])\}\nonumber\\
&=0,
\end{align}
which implies the first equation of (5.7). The rest of (5.7) are clear. 
\par
By substituting (5.7) into (5.5), it turns out that  
\begin{align}
\frac12\Delta\vert d\varphi\vert^2
&=
\sum_{i,j=1}^m\langle (\widetilde{\nabla}_{e_i}
\widetilde{\nabla}_{e_j}d\varphi)(e_j),d\varphi(e_i))\nonumber\\
&-\sum_{i,j=1}^m\langle{}^N\!R(d\varphi(e_i),d\varphi(e_j))d\varphi(e_j),d\varphi(e_i)\rangle\nonumber\\
&+\vert \widetilde{\nabla}d\varphi\vert^2. 
\end{align}
For the first term of (5.11),  using (5.6) and 
$\tau(\varphi)=\sum_{j=1}^m(\widetilde{\nabla}_{e_j}d\varphi)(e_j)$, 
we have 
 \begin{align}
 \sum_{j=1}^m(\widetilde{\nabla}_{e_i}
 \widetilde{\nabla}_{e_j}d\varphi)(e_j)
 &=
 \sum_{j=1}^m\bigg\{
 \overline{\nabla}_{e_i}
 ((\widetilde{\nabla}_{e_j}d\varphi)(e_j))
 -(\widetilde{\nabla}_{e_j}d\varphi)(\nabla_{e_j}e_j)\nonumber\\
 &\qquad\qquad-(\widetilde{\nabla}_{\nabla_{e_i}e_j}d\varphi)(e_j)
 \bigg\}
 \nonumber\\
 &=
 \overline{\nabla}_{e_i}\tau(\varphi)
 -\sum_{j=1}^m(\widetilde{\nabla}_{e_j}d\varphi)(\nabla_{e_i}e_j)\nonumber\\
&
\qquad\qquad -\sum_{j=1}^m(\widetilde{\nabla}_{\nabla_{e_i}e_j}d\varphi)(e_j)
 \end{align}
 Substituting (5.12) into the first term of (5.11), 
 the right hand side of (5.11) is equal to  
 \begin{align}
 \sum_{i,j=1}^m\langle
 \overline{\nabla}_{e_i}\tau(\varphi),d\varphi(e_i)\rangle
&- \sum_{i,j=1}^m
 \{
 \langle(\widetilde{\nabla}_{e_j}d\varphi)(\nabla_{e_i}e_j),d\varphi(e_i)\rangle
 \nonumber\\
 &\qquad\qquad
 +\langle(\widetilde{\nabla}_{\nabla_{e_i}e_j}d\varphi)(e_j),d\varphi(e_i)\rangle
 \}\nonumber\\
 & -\sum_{i,j=1}^m
 \langle {}^N\!R(d\varphi(e_i),d\varphi(e_j))d\varphi(e_j),d\varphi(e_i)\rangle
 \nonumber\\
 &+\vert \widetilde{\nabla}d\varphi\vert^2,
 \end{align}
 which is the desired. 
 \end{pf}
 \vskip0.6cm\par
 \begin{slem} 
 Assume that $(M,g)$ is an 
 $m$-dimensional compact Riemannian manifold and the sectional curvature of $(N,h)$ is bounded from above by a positive constant $C>0$. Then, 
 there exists a positive constant $C'>0$ 
 depending only on $C$ and the geometry of $(M,g)$ such that, 
 for every smooth map $\varphi$ of $M$ into $N$, 
 \begin{align}
 \frac12\Delta\vert d\varphi\vert^2+C'\vert d\varphi\vert^2
 +C'\vert d\varphi\vert^4
 \geq
 \sum_{i=1}^m\langle\overline{\nabla}_{e_i}\tau(\varphi),d\varphi(e_i)\rangle. 
 \end{align} 
 \end{slem}
 \begin{pf} \quad Since the sectional curvature of $(N,h)$ 
 is bounded above by a positive constant $C>0$, 
 \begin{align}
 \langle {}^N\!R(d\varphi(e_i),d\varphi(e_j))d\varphi(e_j),d\varphi(e_i)\rangle
 &\leq 
 C\{
 \vert d\varphi(e_i)\vert^2\vert d\varphi(e_j)\vert^2
 \nonumber\\
 &\qquad
 -\langle d\varphi(e_i),d\varphi(e_j)\rangle^2\}\nonumber\\
 &\leq C\vert d\varphi\vert^4, 
 \end{align}
 and (5.4) in Sublemma 5.2, we have 
 \begin{align}
 \frac12\Delta \vert d\varphi\vert^2
 +C\vert d\varphi\vert^4
 &\geq 
 \vert\widetilde{\nabla}d\varphi\vert^2+
 \sum_{i=1}^m
 \langle\overline{\nabla}_{e_i}\tau(\varphi),d\varphi(e_i)\rangle\nonumber\\
&\quad-\sum_{i,j=1}^m\langle
(\widetilde{\nabla}_{\nabla_{e_i}e_j}d\varphi)(e_j),d\varphi(e_i))
\nonumber\\
&\quad-\sum_{i,j=1}^m
\langle (\widetilde{\nabla}_{e_j}d\varphi)(\nabla_{e_i}e_j),d\varphi(e_i)
\rangle. 
 \end{align}
 Since $(M,g)$ is compact, 
 for some positive constant $C'>0$, the third term of (5.16) 
 is bigger than or equal to 
 \begin{align}
-C'\vert \widetilde{\nabla}d\varphi\vert \,\vert d\varphi\vert
 &\geq 
 C'\varepsilon^2\,\vert \widetilde{\nabla}d\varphi\vert^2-\frac{C'}{4}\frac{1}{\varepsilon^2}\,\vert d\varphi\vert^2 \quad (\forall \,\,\varepsilon>0)
 \nonumber\\
 &
 \geq
 -\frac12\vert \widetilde{\nabla}d\varphi\vert^2
 -\frac{{C'}^2}{2}\,\vert d\varphi\vert^2,
 \end{align}
 by taking 
 $\varepsilon=\frac{1}{\sqrt{2C'}}$. 
 The fourth term of (5.16) is also bigger than or equal to 
 \begin{align}
 -C'\vert \widetilde{\nabla}d\varphi\vert\,\vert d\varphi\vert
 \geq 
 -\frac12\vert\widetilde{\nabla}d\varphi\vert^2
 -\frac{{C'}^2}{2}\,\vert d\varphi\vert^2. 
 \end{align}
 Therefore, the right hand side of (5.16) 
 is bigger than or equal to 
 \begin{align}
 \vert\widetilde{\nabla} d\varphi\vert^2
 &+\sum_{i=1}^m\langle
 \overline{\nabla}_{e_i}\tau(\varphi),d\varphi(e_i)\rangle
 -\frac12\vert\widetilde{\nabla}d\varphi\vert^2
 -\frac{{C'}^2}{2}\,\vert d\varphi\vert^2
 \nonumber\\
 &
 \qquad\qquad\qquad-\frac12\vert\widetilde{\nabla}d\varphi\vert^2
 -\frac{{C'}^2}{2}\vert d\varphi\vert^2\nonumber\\
 &=
 \sum_{i=1}^m
 \langle \overline{\nabla}_{e_i}\tau(\varphi),d\varphi(e_i)\rangle-{C'}^2\vert d\varphi\vert^2, 
 \end{align}
 so that we have 
 \begin{align}
 \frac12\Delta\vert d\varphi\vert^2
 +{C'}^2\,\vert d\varphi\vert^2+C\vert d\varphi\vert^4
 \geq
 \sum_{i=1}^m
 \langle \widetilde{\nabla}_{e_i}\tau(\varphi), 
 d\varphi(e_i)\rangle, 
 \end{align}
 which is (5.14). 
 \end{pf}
\vskip0.6cm\par
\subsection{Three lemmas for $L^{1,q}$-estimates}
Let $r_1$ and $r_2$ be arbitrarily fixed positive real numbers such as
$0<r_1<r_2<1$, and for 
a fixed $x_0\in M$, let us consider a cut off 
function 
$\eta$ satisfying that 
\begin{equation}
\left\{
\begin{aligned}
0&\leq\eta(x)\leq1 \qquad(x\in M),\nonumber\\
&\eta(x)=1\qquad\quad\,\,(x\in B_{r_1}(x_0)),\nonumber\\
&\eta(x)=0\qquad\quad\,\, (x\not\in B_{r_2}(x_0)),\nonumber\\
&\vert\nabla\eta\vert\leq \frac{2}{r_2-r_1}\quad\,\,(\text{everywhere on}\,\,M). 
\end{aligned}
\right.
\end{equation}
Let us also use the notation 
$\overline{\gamma}:=\frac{m}{m-2}$ 
where $m=\dim M$. 
\begin{lem}
There exist constants $D>0$  
depending only on $C>0$ and 
the geometry of $(M,g)$ 
such that, 
for every smooth map $\varphi$ 
of $M$ into $N$, and every $p\geq 3$, 
 \begin{align}
 \left\{
 \int_M(\vert d\varphi\vert^p\,\eta^2)^{\overline{\gamma}}\,v_g
 \right\}^{1/\overline{\gamma}}
 &\leq 
 \frac{p^2D}{(r_2-r_1)^2}
 \int_{B_{r_2}(x_0)}\vert d\varphi\vert^p\,v_g
 \nonumber\\
 &\quad+p\,D\,\int_M
 \vert d\varphi\vert^{p+2}\,\eta^2\,v_g
 \nonumber\\
 &\quad+
 \frac{p\,D}{(r_2-r_1)^2}
 \int_{B_{r_2}(x_0)}\vert\tau(\varphi)\vert^p\,v_g. 
 \end{align}
\end{lem}
\begin{pf}The proof is divided into nine steps.\par 
({\it The first step})\quad 
Multiplying by
$\vert d\varphi\vert^{p-2}\,\eta^2$, (5.14) in 
Sublemma 5.3, and integrating it over $M$, 
we have 
\begin{align}
\int_M\frac12\vert d\varphi\vert^{p-2}\,\eta^2\,&\Delta\vert d\varphi\vert^2\,v_g
+C'\int_M\vert d\varphi\vert^p\,\eta^2\,v_g+
C'\int_M\vert d\varphi\vert^{p+2}\,\eta^2\,v_g\nonumber\\
&\geq 
\int_M\vert d\varphi\vert^{p-2}\,\eta^2\sum_{i=1}^m
\langle 
\overline{\nabla}_{e_i}\tau(\varphi),d\varphi(e_i)\rangle\,v_g. 
\end{align}
\par
({\it The second step})\quad 
For the first term of the left hand side of (5.22), 
we will show that 
\begin{align}
\int_M\frac12\vert d\varphi\vert^{p-2}\,\eta^2\,
\Delta\vert d\varphi\vert^2\,v_g
&\leq 
-\frac{3(p-2)}{p^2}\int_M\vert 
\nabla(\,\vert d\varphi\vert^{p/2}\,)\vert^2\,\eta^2\,v_g
\nonumber\\
&+\frac{16}{p-2}\frac{1}{(r_2-r_1)^2}
\int_{B_{r_2}(x_0)}\vert d\varphi\vert^p\,v_g. 
\end{align}
Indeed, 
integrating the first term by part, 
and using 
$\nabla(\,\vert d\varphi\vert^{p-2}\,\eta^2\,)
=(p-2) \vert d\varphi\vert^{p-3}\eta^2\nabla
(\,\vert d\varphi\vert\,)+2\vert d\varphi\vert^{p-2}\eta\nabla\eta$ 
and 
$\nabla(\,\vert d\varphi\vert^2\,)=2\vert d\varphi\vert\,\nabla(\,\vert d\varphi\vert\,)$, 
we have 
 \begin{align}
 \int_M\frac12
 \vert d\varphi\vert^{p-2}\,\eta^2\,
\Delta\vert d\varphi\vert^2\,v_g
&=-
\int_M\frac12\langle \nabla(\,\vert d\varphi\vert^{p-2}\,\eta^2\,), \nabla(\,\vert d\varphi\vert^2\,)\rangle\,v_g
\nonumber\\
&=
-(p-2)\int_M\vert d\varphi\vert^{p-2}\,\vert\nabla(\,\vert d\varphi\vert\,)\vert^2\,\eta^2\,v_g\nonumber\\
&\quad-2\int_M\langle\nabla\eta,\nabla(\,\vert d\varphi\vert\,)\rangle\,\vert d\varphi\vert^{p-1}\,\eta\,v_g. 
 \end{align}
 Here, substituting both
 \begin{align*}
 \vert d\varphi\vert^{p-2}\,\vert\nabla(\,\vert d\varphi\vert\,)\vert^2
 &=\frac{4}{p^2}
  \vert \nabla(\,\vert d\varphi\vert^{p/2}\,)\vert^2,\\
 \langle\nabla\eta,\nabla(\,\vert d\varphi\vert\,)
 \rangle\,\vert d\varphi\vert^{p-1}\,
 \eta
 &=\frac{2}{p}\,\langle
 \vert d\varphi\vert^{p/2}\,\nabla\eta,
 \eta\,\nabla(\,\vert d\varphi\vert^{p/2}\,)
 \rangle, 
 \end{align*}
 into the right hand side of (5.24),  
 (5.24) turns out that 
 \begin{equation}
 -\frac{4(p-2)}{p^2}\int_M
 \vert\nabla(\,\vert d\varphi\vert^{p/2}\,)\vert^2\,\eta^2\,v_g
 -\frac{4}{p}\int_M
 \langle
  \vert d\varphi\vert^{p/2}\,\nabla\eta,
 \eta\,\nabla(\,\vert d\varphi\vert^{p/2}\,)
 \rangle\,v_g. 
 \end{equation}
 Here, for $A:=\sqrt{\frac{p}{p-2}}\vert d\varphi\vert^{p/2}\,\nabla\eta$, 
 $B:=\sqrt{\frac{p-2}{p}}\eta\,\nabla(\,\vert d\varphi\vert^{p/2}\,)$, 
 it holds that 
 \begin{align}
 -\langle
 &\vert d\varphi\vert^{p/2}\nabla\eta,\eta\nabla(\,\vert d\varphi\vert^{p/2}\,)
 \rangle\nonumber\\
&= -\langle A, B\rangle \nonumber\\
&\leq\vert A\vert^2+\frac14\vert B\vert^2\nonumber\\
&=\frac{p}{p-2}\vert d\varphi\vert^p\,\vert\nabla\vert^2+\frac14\frac{p-2}{p}\eta^2\,\vert\nabla(\,\vert d\varphi\vert\,)\vert^2, 
 \end{align}
 so that, for the second term of (5.25),  
 \begin{align}
 -\frac{4}{p}&\int_M
 \langle
 \vert d\varphi\vert^{p/2}\nabla\eta,\eta\,\nabla(\,\vert d\varphi\vert^{p/2}\,)
 \rangle\,v_g\nonumber\\
 &\leq
 \frac{4}{p}\frac{p}{p-2}\int_M\vert d\varphi\vert^p\,\vert\nabla\eta\vert^2\,v_g
 +\frac{4}{p}\frac{p-2}{4p}\int_M
 \eta^2\,\vert\nabla(\,\vert d\varphi\vert^{p/2}\,)\vert^2\,v_g\nonumber\\
 &=
 \frac{4}{p-2}\int_M\vert d\varphi\vert^p\,\vert\nabla\eta\vert^2\,v_g
 +\frac{p-2}{p^2}\int_M
 \eta^2\,\vert\nabla(\,\vert d\varphi\vert^{p/2}\,)\vert^2\,v_g. 
 \end{align}
 Therefore, we have 
 \begin{align}
 \int_M&\frac12\vert d\varphi\vert^{p-2}\,\eta^2\,\Delta(\,\vert d\varphi\vert^2\,)\,v_g
 \nonumber\\
 &\leq\left\{-\frac{4(p-2)}{p^2}
 +\frac{p-2}{p^2}\right\}
 \int_M
 \vert\nabla(\,\vert d\varphi\vert^{p/2}\,)\vert^2\,\eta^2\,v_g\nonumber\\
 &\qquad\qquad+\frac{4}{p-2}\int_M
  \vert d\varphi\vert^{p/2}
  \,\vert\nabla\eta\vert^2
  \,v_g\nonumber\\
  &\leq
  -\frac{3(p-2)}{p^2}\int_M\vert 
\nabla\,\vert d\varphi\vert^{p/2}\,\vert^2\,\eta^2\,v_g
\nonumber\\
&\qquad\qquad
+\frac{16}{p-2}\frac{1}{(r_2-r_1)^2}
\int_{B_{r_2}(x_0)}\vert d\varphi\vert^p\,v_g, 
 \end{align}
 which implies (5.23). 
 \par
 ({\it The third step})\quad 
 For the right hand side of (5.22), we have 
 \begin{align}
 \int_M
 &\vert d\varphi\vert^{p-2}\,\eta^2\sum_{i=1}^m
\langle 
\overline{\nabla}_{e_i}\tau(\varphi),d\varphi(e_i)\rangle\,v_g\nonumber\\
&=\int_M
\sum_{i=1}^m
\langle
\overline{\nabla}_{e_i}\tau(\varphi), 
\vert d\varphi\vert^{p-2}\,\eta^2\,d\varphi(e_i)\rangle\,v_g\nonumber\\
&=-\int_M
\langle\tau(\varphi),
\sum_{i=1}^m
\widetilde{\nabla}_{e_i}
(\,\vert d\varphi\vert^{p-2}\,\eta^2\,d\varphi\,)(e_i)
\rangle\,v_g\nonumber\\
&=-\int_M\langle\tau(\varphi), \sum_{i=1}^m 
\{e_i(\,\vert d\varphi\vert^{p-2}\,)\,\eta^2\,d\varphi(e_i)
\}\rangle\,v_g\nonumber\\
&\quad-\int_M\langle\tau(\varphi), 
\vert d\varphi\vert^{p-2}\sum_{i=1}^m
\{e_i(\,\eta^2\,)\,d\varphi(e_i)
\}\rangle\,v_g\nonumber\\
&\quad-\int_M\langle\tau(\varphi), 
\vert d\varphi\vert^{p-2}\,\eta^2\,\tau(\varphi)\rangle\,v_g\nonumber\\
&\geq
-\int_M\vert\tau(\varphi)\vert\,\vert\nabla(\,\vert d\varphi\vert^{p-2}\,)\vert\,\,\eta^2\,
\vert d\varphi\vert\,v_g\nonumber\\
&\quad -\int_M 
\vert\tau(\varphi)\vert\,\vert d\varphi\vert^{p-1}\,\vert\nabla(\,\eta^2\,)\vert\,
v_g
-\int_M\vert\tau(\varphi)\vert^2\,\vert d\varphi\vert^{p-2}\,\eta^2\,v_g. 
 \end{align}
 Here, for the first term of the last right hand side of (5.29), we have that 
 \begin{align}
 \vert\tau(\varphi)\vert\,&\vert\nabla(\,\vert d\varphi\vert^{p-2}\,)\vert\,\,\eta^2\,
\vert d\varphi\vert
=(p-2) \vert\tau(\varphi)\vert\,\vert d\varphi\vert^{p-2}\,\vert
\nabla(\,\vert d\varphi\vert\,)\vert\,\eta^2
\nonumber\\
&=\frac{2(p-2)}{p}
\left(
\sqrt{\frac{1}{p}}\,\vert\nabla(\,\vert d\varphi\vert^{p/2}\,)\vert\,\eta
\right)
\left(
\sqrt{p}\,\vert d\varphi\vert^{p/2-1}\,\eta\,\vert\tau(\varphi)\vert
\right)\nonumber\\
&\leq 
\frac{2(p-2)}{p}
\bigg\{
\frac12\frac{1}{p}\,\vert\,\nabla(\,\vert d\varphi\vert^{p/2}\,)\vert^2\,\eta^2+
\frac12\,p\,\vert d\varphi\vert^{p-2}\,\eta^2\,\vert\tau(\varphi)\vert^2
\bigg\}
\nonumber\\
&=
\frac{p-2}{p^2}\vert\nabla(\,\vert d\varphi\vert^{p/2}\,)\vert^2\,\eta^2
+(p-2)\,\vert d\varphi\vert^{p-2}\,\eta^2\,\vert\tau(\varphi)\vert^2. 
  \end{align}
  Substituting (5.30) into the first term of 
  the right hand side of (5.29), 
  (5.29) turns out that  
   \begin{align}
   \int_M
 &\vert d\varphi\vert^{p-2}\,\eta^2\sum_{i=1}^m
\langle 
\overline{\nabla}_{e_i}\tau(\varphi),d\varphi(e_i)\rangle\,v_g\nonumber\\
&\geq 
-\frac{p-2}{p^2}\int_M 
\vert \,\nabla(\,\vert d\varphi\vert^{p/2}\,)\vert^2\,\eta^2\,v_g
-(p-2)\int_M\vert d\varphi\vert^{p-2}\,\eta^2\,\vert \tau(\varphi)\vert^2\,v_g\nonumber\\
&\quad
-\int_M\vert\tau(\varphi)\vert\,\vert d\varphi\vert^{p-1}\,2\eta\,\vert\nabla\eta\vert\,v_g
-\int_M\vert\tau(\varphi)\vert^2\,\vert d\varphi\vert^{p-2}\,\eta^2\,v_g\nonumber\\
&=-\frac{p-2}{p^2}
\int_M\vert \,\nabla(\,\vert d\varphi\vert^{p/2}\,)\vert^2\,\eta^2\,v_g
-(p-1)\int_M\vert d\varphi\vert^{p-2}\,\eta^2\,\vert \tau(\varphi)\vert^2\,v_g\nonumber\\
&\quad 
-2\int_M\vert d\varphi\vert^{p-1}\,\eta\,\vert\nabla\eta\vert\,\vert\tau(\varphi)\vert\,v_g.
  \end{align}
  \par
  ({\it The fourth step})\quad
  For the right hand side of (5.31), we have
  \begin{align}
  \int_M\vert d\varphi\vert^{p-2}\,\eta^2\,\vert\tau(\varphi)\vert^2\,v_g
  &\leq\frac{p-2}{p}\int_{B_{r_2}(x_0)}\vert d\varphi\vert^p\,v_g\nonumber\\
  &\quad
  +\frac{2}{p}\int_{B_{r_2}(x_0)}\vert\tau(\varphi)\vert^p\,v_g,\\
  \int_M\vert d\varphi\vert^{p-1}\,\eta\,\vert\nabla\eta\vert\,\vert\tau(\varphi)\vert\,v_g
  &\leq \frac{p-1}{p}\frac{2}{r_2-r_1}\int_{B_{r_2}(x_0)}\vert d\varphi\vert^p\,v_g\nonumber\\
  &\quad+\frac{1}{p}\,\frac{2}{r_2-r_1}\int_{B_{r_2}(x_0)}\vert\tau(\varphi)\vert^p\,v_g.
  \end{align}
  \par
  Indeed, for (5.32), we have, 
  by using H\"older's inequality
  because of  $\frac{p-2}{p}+\frac{2}{p}=1$, 
  \begin{align}
  \int_M\vert d\varphi\vert^{p-2}\,\eta^2\,\vert\tau(\varphi)\vert^2\,v_g
  &\leq\int_{B_{r_2}(x_0)}\vert d\varphi\vert^{p-2}\,\vert\tau(\varphi)\vert^2\,v_g
  \nonumber\\
  &\leq
  \left(
  \int_{B_{r_2}(x_0)}\vert d\varphi\vert^p\,v_g
  \right)^{(p-2)/p}
  \nonumber\\
  &\qquad\times
  \,\left(
  \int_{B{r_2}(x_0)}
  \vert\tau(\varphi)\vert^p\,v_g
  \right)^{2/p}
  \nonumber\\
  &\leq
  \frac{p-2}{p}\int_{B_{r_2}(x_0)}
  \vert d\varphi\vert^p\,v_g+\frac{2}{p}
  \int_{B_{r_2}(x_0)}\vert\tau(\varphi)\vert^p\,v_g
  \nonumber
  \end{align}
  which is (5.32). 
  In the last inequality, we used the inequality: 
  $AB\leq \frac{p-2}{p}\,A^{p/(p-2)}+\frac{2}{p}\,B^{p/2}$ for every positive real numbers $A$ and $B$. 
  \par
  For (5.33), we have again using H\"older's inequality since 
  $\frac{p-1}{p}+\frac{1}{p}=1$, 
  \begin{align}
  \int_M&\vert d\varphi\vert^{p-1}\,\eta\,\vert\nabla\eta\vert\,\vert\tau(\varphi)\vert\,v_g
\leq 
\frac{2}{r_2-r_1}\int_{B_{r_2}(x_0)}
\vert d\varphi\vert^{p-1}\,\vert\tau(\varphi)\vert\,v_g
\nonumber\\
&\leq
\frac{2}{r_2-r_1}
\left(
\int_{B_{r_2}(x_0)}\vert d\varphi\vert^p\,v_g
\right)^{(p-1)/p}
\,
\left(
\int_{B_{r_2}(x_0)}\vert\tau(\varphi)\vert^p
\right)^{1/p}\nonumber\\
&\leq
\frac{p-1}{p}\frac{2}{r_2-r_1}\int_{B_{r_2}(x_0)}
\vert d\varphi\vert^p\,v_g
+\frac{1}{p}\,\,\frac{2}{r_2-r_1}
\int_{B_{r_2}(x_0)}\vert\tau(\varphi)\vert^p\,v_g
\nonumber
  \end{align}
  which is (5.33). In the last inequality, 
  we used again the inequality: 
  $AB\leq \frac{p-1}{p}\,A^{p/(p-1)}+\frac{1}{p}\,B^p$. 
  \par
  ({\it The fifth step})\quad 
  By substituting 
  (5.32) and (5.33) into (5.31), we have, 
  for $p\geq 3$, 
  \begin{align}
   &\int_M
 \vert d\varphi\vert^{p-2}\,\eta^2\sum_{i=1}^m
\langle 
\overline{\nabla}_{e_i}\tau(\varphi),d\varphi(e_i)\rangle\,v_g
\geq \nonumber\\
&-\frac{p-2}{p^2}\int_M\vert\nabla(\,\vert d\varphi\vert^{p/2}\,)\vert^2\,\eta^2\,v_g\nonumber\\
&-\big\{\frac{(p-1)(p-2)}{p}+\frac{p-1}{p}\,\frac{4}{r_2-r_1}\big\}
\int_{B_{r_2}(x_0)}\vert d\varphi\vert^p\,v_g
\nonumber\\
&-\big\{
\frac{2(p-1)}{p}+\frac{1}{p}\,\frac{4}{r_2-r_1}
\big\}
\int_{B_{r_2}(x_0)}\vert\tau(\varphi)\vert^p\,v_g. 
  \end{align}
  \par
  ({\it The sixth step})\quad 
  Together (5.22), (5.23) and (5.34), we obtain 
  \begin{align}
 & \frac{2(p-2)}{p^2}
  \int_M
  \vert\nabla(\,\vert d \varphi\vert^{p/2}\,)\vert^2\,\eta^2\,v_g
  \leq \nonumber\\
  &\big\{C'+\frac{p-1}{p}(p-2)
  +\frac{p-1}{p}\,\frac{4}{r_2-r_1}
  +\frac{16}{p-2}\,\frac{1}{(r_2-r_1)^2}\big\}
  \int_{B_{r_2}(x_0)}\vert d\varphi\vert^p\,v_g
  \nonumber\\
 &\quad+C'\int_M\vert d\varphi\vert^{p+2}\,\eta^2\,v_g
  \nonumber\\
  &\quad+\big\{\frac{2(p-1)}{p}+\frac{1}{p}\,\frac{4}{r_2-r_1}\big\}
  \int_{B_{r_2}(x_0)}\vert\tau(\varphi)\vert^p\,v_g,
  \end{align} 
  so that we have 
  \begin{align}
  \int_M\vert\nabla(\,\vert d\varphi\vert^{p/2}\,)\vert^2\,\eta^2\,v_g
  &\leq A_1
  \int_{B_{r_2}(x_0)}\vert d\varphi\vert^p\,v_g
  +A_2\int_M
  \vert d\varphi\vert^{p+2}\,\eta^2\,v_g
  \nonumber\\
  &\quad +A_3
  \int_{B_{r_2}(x_0)}\vert\tau(\varphi)\vert^p\,v_g,
  \end{align}
  where 
  \begin{align}
  A_1&=
  \frac{p^2}{2(p-2)}\,\big\{
  C'+\frac{p-1}{p}(p-2)
  +\frac{p-1}{p}\,\frac{4}{r_2-r_1}+\frac{16}{p-1}\,\frac{1}{(r_2-r_1)^2}
 \big \}\nonumber\\
  &\leq\frac{p^2(C'+4)}{(r_2-r_1)^2},\\
  A_2&=\frac{p^2}{2(p-2)}\,C'\leq p\,C',\\
  A_3&=
  \frac{p^2}{2(p-2)}\,\big\{
  \frac{2(p-1)}{p}+\frac{4}{p}\,\frac{1}{r_2-r_1}
  \big\}\leq 
  \frac{4p}{(r_2-r_1)^2}, 
  \end{align}
  since $p\geq 3$ and $0<r_2-r_1<1$. 
  \par
  ({\it The seventh step})\quad 
  On the other hand, since
  \begin{align*}
  \vert\nabla(\,\vert d\varphi\vert^{p/2}\,\eta\,)\vert
  &=\vert
  \eta\,\nabla(\,\vert d\varphi\vert^{p/2}\,)+
  \vert d\varphi\vert^{p/2}\,\nabla\eta\vert^2\\
  &\leq
  2\,\vert\nabla(\,\vert d\varphi\vert^{p/2}\,)\vert^2\,
  \eta^2
  +2\,\vert d\varphi\vert^p\,\vert\nabla\eta\vert^2,
 \end{align*}
 we have 
 \begin{align}
 \int_M
 \vert\nabla(\vert d\varphi\vert^{p/2}\,\eta\,)\vert^2\,
 v_g
& \leq
 2\int_M
 \vert\nabla(\,\vert d\varphi\vert^{p/2}\,)\vert^2\,\eta^2\,v_g\nonumber\\
 &\quad+\frac{8}{(r_2-r_1)^2}
 \int_{B_{r_2}(x_0)}\vert d\varphi\vert^p\,v_g. 
 \end{align}
 \par
 ({\it The eighth step})\quad 
 By Sobolev's inequality 
 for a function $\vert d\varphi\vert^{p/2}\,\eta$ 
 on $(M,g)$, and $\overline{\gamma}=\frac{m}{m-2}$, 
 we have 
  \begin{align*}
 \bigg(
  \int_M
  (\,\vert d\varphi\vert^p\,\eta^2\,)
  ^{\overline{\gamma}}\,v_g
  \bigg)
  ^{1/(2\overline{\gamma})}
  &=\bigg(
  \int_M(\,\vert d\varphi\vert^{p/2}\,\eta\,)
  ^{2m/(m-2)}\,v_g
  \bigg)^{(m-2)/(2m)}\nonumber\\
  &\leq
  C''\bigg(
  \int_M
  \vert\nabla(\,\vert d\varphi\vert^{p/2}\,\eta\,)\vert^2\,v_g
  \bigg)^{1/2},
  \end{align*}
  so that we have 
  \begin{equation}
   \bigg(
  \int_M
  (\,\vert d\varphi\vert^p\,\eta^2\,)
  ^{\overline{\gamma}}\,v_g
  \bigg)
  ^{1/\overline{\gamma}}
  \leq
  C''
  \int_M
  \vert\nabla(\,\vert d\varphi\vert^{p/2}\,\eta\,)\vert^2\,v_g. 
  \end{equation}
  \par
  ({\it The ninth step})\quad 
  Together with (5.36)\,$\sim$ (5.41), 
  we have 
  \begin{align*}
    \bigg(
  \int_M
  &(\,\vert d\varphi\vert^p\,\eta^2\,)
  ^{\overline{\gamma}}\,v_g
  \bigg)
  ^{1/\overline{\gamma}}\\
  &\leq
  C''\int_M
  \vert\nabla(\,\vert d\varphi\vert^{p/2}\,\eta\,)\vert^2\,v_g\\
  &
  \leq
  2\,C''\int_M\vert\nabla
  (\,\vert d\varphi\vert^{p/2}\,)\vert^2\,\eta^2\,v_g
 +\frac{8C''}{(r_2-r_1)^2}
  \int_{B_{r_2}(x_0)}\vert d\varphi\vert^p\,v_g
  \\
 & 
 \leq\{
  2C''A_1
  +\frac{8C''}{(r_2-r_1)^2}
  \}\,
  \int_{B_{r_2}(x_0)}\vert d\varphi\vert^p\,v_g
  \\
  &
  \quad+2C''A_2\int_M\vert d\varphi\vert^{p+2}\,\eta^2\,v_g+
  2C''A_3
  \int_{B_{r_2}(x_0)}
  \vert\tau(\varphi)\vert^p\,v_g\\
  &\leq
  \frac{p^2\,(2C''(C'+8))}{(r_2-r_1)^2}\int_{B_{r_2}(x_0)}\vert d\varphi\vert^p\,v_g\\
  &\quad+p\,(2C'C'')\int_M\vert d\varphi\vert^{p+2}\,\eta^2\,v_g
  +\frac{p\,(8C'')}{(r_2-r_1)^2}\int_{B_{r_2}(x_0)}\vert\tau(\varphi)\vert^p\,v_g,
  \end{align*}
  which yields the desired inequality (5.21) by 
  taking $D=\max\{2C''(C'+8), 2C'C'', 8C''\}$. 
  We have Lemma 5.4. 
\end{pf}
\vskip0.6cm\par
Next, we will show
\begin{lem}
Assume that $m=\dim M\geq 3$. Then, 
there exist positive real numbers 
$\varepsilon_i$ $(i=1,2)$ such that, 
if $\int_{B_{\rho_2}(x_0)}\vert d\varphi\vert^m\,v_g\leq \varepsilon_1$ and 
$\int_{B_{\rho_2}(x_0)}\vert\tau(\varphi)\vert^m\,v_g\leq \varepsilon_2$,  
then 
\begin{equation}
\left(
\int_{B_{\frac{\rho_1+\rho_2}{2}}(x_0)}
\vert d\varphi\vert^{m\,\overline{\gamma}}\,v_g
\right)^{1/\overline{\gamma}}
\leq \frac{D'}{(\rho_2-\rho_1)^2},
\end{equation}
for all $\rho_i$ $(i=1,2)$ with 
$0<\rho_1<\rho_2<1$. Here,  $D'=\frac{4mD(m\varepsilon_1+\varepsilon_2)}{1-mD\varepsilon_1{}^{2/m}}$.  
\end{lem}
\begin{pf}
We apply Lemma 5.4 for 
$r_1=\frac{\rho_1+\rho_2}{2}$, $r_2=\rho_2$, 
and $p=m=\dim M\geq 3$. Notice that $r_2-r_1=\frac12(\rho_2-\rho_1)$. 
Then, we have 
\begin{align}
\bigg\{
&\int_M
(\,\vert d\varphi\vert^m\,\eta^2\,)
^{\overline{\gamma}}\,v_g
\bigg\}
^{1/\overline{\gamma}}
\nonumber\\
&\leq
\frac{4m^2D}{(\rho_2-\rho_1)^2}
\int_{B_{\rho_2}(x_0)}
\vert d\varphi\vert^m\,v_g
+mD\,\int_M
\vert d\varphi\vert^{m+2}\,\eta^2\,v_g
\nonumber\\
&\quad 
+\frac{4mD}{(\rho_2-\rho_1)^2}
\int_{B_{\rho_2}(x_0)}\vert\tau(\varphi)\vert^m\,v_g
\nonumber\\
&\leq 
\frac{4m^2D\varepsilon_1}{(\rho_2-\rho_1)^2}
+\frac{4mD\varepsilon_2}{(\rho_2-\rho_1)^2}
+mD\int_M\vert d\varphi\vert^{m+2}\,\eta^2\,v_g. 
\end{align}
On the other hand, since 
$\frac{1}{\overline{\gamma}}+\frac{2}{m}=1$, by H\"older's inequality, 
\begin{align}
\int_M&\vert d\varphi\vert^{m+2}\,\eta^2\,v_g
=\int_{B_{\rho_2}(x_0)}(\,\vert d\varphi\vert^m\,\eta^2\,)\,\vert d\varphi\vert^2\,v_g\nonumber\\
&\leq
\bigg\{
\int_{B_{\rho_2}(x_0)}
(\,\vert d\varphi\vert^m\,\eta^2\,)^{\overline{\gamma}}\,v_g
\bigg\}^{1/\overline{\gamma}}\,\,
\bigg\{
\int_{B_{\rho_2}(x_0)}
(\,\vert d\varphi\vert^2\,)^{m/2}
\bigg\}^{2/m}\nonumber\\
&\leq \varepsilon_1{}^{2/m}
\,\bigg\{
\int_M(\,\vert d\varphi\vert^m\,\eta^2\,)^{\overline{\gamma}}\,v_g
\bigg\}^{1/\overline{\gamma}}.
\end{align}
Therefore, we obtain by (5.43) and (5.44), 
\begin{align*}
\bigg\{
\int_M(\,\vert d\varphi\vert^m\,\eta^2\,)^{\overline{\gamma}}\,v_g
\bigg\}^{1/\overline{\gamma}}
&\leq
\frac{4mD(m\,\varepsilon_1+\varepsilon_2)}{(\rho_2-\rho_1)^2}\nonumber\\
&\quad+mD\,\epsilon_1{}^{2/m}
\bigg\{
\int_M(\,\vert d\varphi\vert^m\,\eta^2\,)^{\overline{\gamma}}\,v_g
\bigg\}^{1/\overline{\gamma}}, 
\end{align*}
so that 
\begin{align*}
(1-mD\varepsilon_1{}^{2/m})
\,\bigg\{
\int_M(\,\vert d\varphi\vert^m\,\eta^2\,)^{\overline{\gamma}}\,v_g
\bigg\}^{1/\overline{\gamma}}
\leq
\frac{4mD(m\,\varepsilon_1+\varepsilon_2)}{(\rho_2-\rho_1)^2}. 
\end{align*}
Thus, if we choose a positive number 
$\varepsilon_1$ 
in such a way that $\varepsilon_1<(\frac{1}{mD})^{m/2}$, then 
we obtain 
\begin{equation}
\bigg\{
\int_M(\,\vert d\varphi\vert^m\,\eta^2\,)^{\overline{\gamma}}\,v_g
\bigg\}^{1/\overline{\gamma}}
\leq 
\frac{4mD(m\,\varepsilon_1+\varepsilon_2)}{1-mD\varepsilon_1{}^{2/m}}\,
\frac{1}
{(\rho_2-\rho_1)^2}.
\end{equation}
Now, since 
$r_1=\frac{\rho_1+\rho_2}{2}$ and 
$\eta=1$ on $B_{r_1}(x_0)$, we have 
\begin{align}
\bigg\{
\int_{B_{\frac{\rho_1+\rho_2}{2}}(x_0)}
\vert d\varphi\vert^{m\,\overline{\gamma}}\,v_g
\bigg\}^{1/\overline{\gamma}}
&\leq 
\bigg\{
\int_M(\,\vert d\varphi\vert^m\,\eta^2\,)^{\overline{\gamma}}\,v_g
\bigg\}^{1/\overline{\gamma}}\nonumber\\
&\leq
\frac{4mD(m\,\varepsilon_1+\varepsilon_2)}{1-mD\varepsilon_1{}^{2/m}}\,
\frac{1}
{(\rho_2-\rho_1)^2},
\end{align}
which is Lemma 5.5. 
\end{pf}
\vskip0.6cm\par
\begin{lem} Assume that $m=\dim M\geq 3$. 
For every $p\geq 3$, 
and $0<\rho_1<\rho_2<1$, it holds that 
\begin{align}
\bigg\{
\int_{B_{\rho_1}(x_0)}
\vert d\varphi\vert^{p\overline{\gamma}}v_g
\bigg\}^{1/(p\,\overline{\gamma})}
&\leq
\frac{\max(p^{2/p},p^{m/(2p)}){D''}^{1/p}}
{
(\rho_2-\rho_1)^{2/p}}
\bigg\{
\int_{B_{\rho_2}(x_0)}\vert d\varphi\vert^pv_g
\bigg\}^{1/p}\nonumber\\
&\quad +
\frac{p^{1/p}\,{D''}^{1/p}}{(\rho_2-\rho_1)^{1/p}}
\bigg\{
\int_{B_{\rho_2}(x_0)}
\vert\tau(\varphi)\vert^p\,v_g
\bigg\}^{1/p},
\end{align}
where $D''=\max(8(D+\frac{\overline{\gamma}}{m}\,(DD')^{m/2},\,8D)$. 
\end{lem}
\begin{pf}
\quad The proof is divided into three steps. 
\par
({\it The first step})
\quad 
We apply Lemma 5.4 in the case that 
$r_1=\rho_1$ and $r_2=\frac{\rho_1+\rho_2}{2}$ 
which implies that 
$r_2-r_1=\frac12(\rho_2-\rho_1)$.  
Then, we obtain 
\begin{align}
\left\{
 \int_M(\vert d\varphi\vert^p\,\eta^2)
 ^{\overline{\gamma}}\,v_g
 \right\}^{1/\overline{\gamma}}
 &\leq 
 \frac{4p^2D}{(\rho_2-\rho_1)^2}
 \int_{B_{\frac{\rho_1+\rho_2}{2}}(x_0)}\vert d\varphi\vert^p\,v_g
 \nonumber\\
 &\quad+p\,D\,\int_{B_{\frac{\rho_1+\rho_2}{2}}(x_0)}
 \vert d\varphi\vert^{p+2}\,\eta^2\,v_g
 \nonumber\\
&\quad+
 \frac{4p\,D}{(\rho_2-\rho_1)^2}
 \int_{B_{\frac{\rho_1+\rho_2}{2}}(x_0)}\vert\tau(\varphi)\vert^p\,v_g. 
\end{align}
Since $\frac{\rho_1+\rho_2}{2}<\rho_2$, 
the sum of the first and last terms of the 
right hand side of (5.48) is smaller than or equal to 
\begin{equation}
\frac{4p^2D}{(\rho_2-\rho_1)^2}
 \int_{B_{\rho_2}(x_0)}\vert d\varphi\vert^p\,v_g
 +\frac{4p\,D}{(\rho_2-\rho_1)^2}
 \int_{B_{\rho_2}(x_0)}\vert\tau(\varphi)\vert^p\,v_g. 
\end{equation}
\par
({\it The second step})\quad 
We estimate the second term of the right hand side of (5.48) as follows: 
For this, we put 
$a:=\frac{m}{2}\,\overline{\gamma}=\frac{m}{2}\,\frac{m}{m-2}$, 
$b:=\overline{\gamma}^2=(\frac{m}{m-2})^2$, and 
$c:=\frac{m}{2}$. 
Since $\frac{\overline{\gamma}}{b}+\frac{1}{c}
=\frac{1}{\overline{\gamma}}+\frac{1}{c}=
\frac{m-2}{m}+\frac{2}{m}=1$, 
the second term 
is equal to 
\begin{equation}
pD\int_{B_{\frac{\rho_1+\rho_2}{2}}(x_0)}
\vert d\varphi\vert^2\,
(\,\vert d\varphi\vert^p\,\eta^2\,)^{\overline{\gamma}/b}\,
(\,\vert d\varphi\vert^p\,\eta^2\,)^{1/c}
\,v_g,
\end{equation}
which, by H\"older's inequality since  $\frac{1}{a}+\frac{1}{b}+\frac{1}{c}=1$, 
is smaller than or equal to 
\begin{align}
&pD\,
\bigg\{
\int_{B_{\frac{\rho_1+\rho_2}{2}}(x_0)}
(\,\vert d\varphi\vert^2\,)^a\,v_g
\bigg\}^{1/a}\,
\bigg\{
\int_{B_{\frac{\rho_1+\rho_2}{2}}(x_0)}
(\,\vert d\varphi\vert^p\,\eta^2\,)^{\overline{\gamma}}\,v_g
\bigg\}^{1/b}\nonumber\\
&\qquad\times\bigg\{
\int_{B_{\frac{\rho_1+\rho_2}{2}}(x_0)}
\vert d\varphi\vert^p\,\eta^2\,v_g
\bigg\}^{1/c}
\nonumber\\
&\leq pD
\,
\bigg\{
\int_{B_{\frac{\rho_1+\rho_2}{2}}(x_0)}
\vert d\varphi\vert^{m\,\overline{\gamma}}\,v_g
\bigg\}^{2/(m\,\overline{\gamma})}\,
\bigg\{
\int_M
(\,\vert d\varphi\vert^p\,\eta^2\,)^{\overline{\gamma}}\,v_g
\bigg\}^{1/\overline{\gamma}^2}
\nonumber\\
&\qquad\times\bigg\{
\int_{B_{\rho_2}(x_0)}
\vert d\varphi\vert^p\,\eta^2\,v_g
\bigg\}^{2/m}
\nonumber\\
&\leq
\frac{pDD'}{(\rho_2-\rho_1)^{4/m}}
\bigg\{
\int_M
(\,\vert d\varphi\vert^p\,\eta^2\,)^{\overline{\gamma}}\,v_g
\bigg\}^{1/\overline{\gamma}^2}
\,
\bigg\{
\int_{B_{\rho_2}(x_0)}
\vert d\varphi\vert^p\,\eta^2\,v_g
\bigg\}^{2/m}
\end{align}
for the last inequality of which follows 
from Lemma 5.5.  
Because of 
$\frac{1}{\overline{\gamma}}+\frac{2}{m}=1$, 
we can apply to (5.51) the Young's inequality 
$AB\leq \varepsilon\,A^{\overline{\gamma}}+
\frac{2\overline{\gamma}}{m}\,\frac{1}{\varepsilon}\,
B^{m/2}$ for every positive real numbers 
$A$, $B$ and $\varepsilon$, 
the last line of (5.51) is smaller than or equal to 
\begin{align}
&\varepsilon \,
\bigg\{
\int_M
(\,\vert d\varphi\vert^p\,\eta^2\,)^{\overline{\gamma}}\,v_g
\bigg\}^{1/\overline{\gamma}}
+\frac{2\overline{\gamma}}{m}\,\frac{1}{\varepsilon}\,
\,
\bigg[
\frac{pDD'}{(\rho_2-\rho_1)^{4/m}}\,
\bigg\{
\int_{B_{\rho_2}(x_0)}
\vert d\varphi\vert^p\,\eta^2\,v_g
\bigg\}^{2/m}
\bigg]^{m/2}\nonumber\\
&=
\varepsilon \,
\bigg\{
\int_M
(\,\vert d\varphi\vert^p\,\eta^2\,)^{\overline{\gamma}}\,v_g
\bigg\}^{1/\overline{\gamma}}
+\frac{2\overline{\gamma}}{m}\,\frac{1}{\varepsilon}\,
\frac{(pDD')^{m/2}}{(\rho_2-\rho_1)^2}\,
\int_{B_{\rho_2}(x_0)}\vert d\varphi\vert^p\,\eta^2\,v_g. 
\end{align}
\par
({\it The third step})\quad 
By (5.48), (5.49) and (5.52), 
we obtain 
\begin{align}
(1-\varepsilon)\,&
\bigg\{
\int_M
(\,\vert d\varphi\vert^p\,\eta^2\,)^{\overline{\gamma}}\,v_g
\bigg\}^{1/\overline{\gamma}}\nonumber\\
&\leq
\big(4p^2D+\frac{2\overline{\gamma}}{m}\,\frac{1}{\varepsilon}\,(pDD')^{m/2}\big)
\,\frac{1}{(\rho_2-\rho_1)^2}\,\int_{B_{\rho_2}(x_0)}\vert d\varphi\vert^p\,v_g\nonumber\\
&\quad
+\frac{4pD}{(\rho_2-\rho_1)^2}\,\int_{B_{\rho_2}(x_0)}\vert\tau(\varphi)\vert^p\,v_g.
\end{align}
By taking $\varepsilon=\frac12$, and noticing that 
$\eta=1$ on $B_{r_1}(x_0)=B_{\rho_1}(x_0)$, 
we have 
\begin{align}
\bigg\{\int_M
(\,\vert d\varphi\vert^p\,\eta^2\,)^{\overline{\gamma}}\,v_g
\bigg\}^{1/\overline{\gamma}}
&\leq
\frac{\max(p^2,\,\,p^{m/2})\,D''}{(\rho_2-\rho_1)^2} 
\,\int_{B_{\rho_2}(x_0)}\vert d\varphi\vert^p\,v_g\nonumber\\
&\quad+\frac{p\,D''}{(\rho_2-\rho_1)^2}\,
\int_{B_{\rho_2}(x_0)}\vert\tau(\varphi)\vert^p\,v_g,
\end{align}
where 
$D''=\max(8(D+\frac{\overline{\gamma}}{m}\,(DD')^{m/2}), 8D)$. 
Finally, taking $\frac{1}{p}$ th-power of the both hand sides of (5.54), and 
using the inequality $(A+B)^{1/p}\leq
A^{1/p}+B^{1/p}$, we have 
Lemma 5.6. 
\end{pf}
\vskip0.6cm\par
\subsection{Proof of Theorem 5.1} 
Now we are in position to give a proof of Theorem 5.1. The proof is carried out again by the Moser's iteration method and divided into five steps. 
\par
({\it The first step})\quad 
We first choose and fix a small positive real number 
$0<\rho^{\ast}<1$ such that 
Vol $(B_{\rho^{\ast}}(x_0))<1$.  Then, for all 
$0<\rho_1<\rho_2<\rho^{\ast}$, we have 
\begin{align}
\bigg\{
\int_{B_{\rho_1}(x_0)}\vert\tau(\varphi)\vert^{p\,\overline{\gamma}}\,v_g
\bigg\}^{1/(p\,\overline{\gamma})}
&\leq \sup_{B_{\rho_1}(x_0)}\vert\tau(\varphi)\vert
\,\,(\text{Vol}(B_{\rho_1}(x_0))^{1/(p\,\overline{\gamma})}\nonumber\\
&\leq \sup_{B_{\rho_1}(x_0)}\vert\tau(\varphi)\vert
\nonumber\\
&\leq
\frac{C'}{(2\,\rho_1)^{m/2}}\,E_2(\varphi)=:D_3<\infty,
\end{align}
by Proposition 4.3. 
\par
({\it The second step})\quad 
We put 
\begin{equation}
\Phi(p,\rho):=
\bigg\{
\int_{B_{\rho}(x_0)}
\vert d\varphi\vert^p\,v_g
\bigg\}^{1/p}
+
\bigg\{
\int_{B_{\rho}(x_0)}
\vert\tau(\varphi)\vert^p\,v_g
\bigg\}^{1/p}
+1.
\end{equation}
Then we want to show that 
\begin{equation}
\Phi(p\,\overline{\gamma},\rho_1)
\leq 
\frac{\max(\,p^{2/p},\,\,p^{m/(2p)}\,)\,D_4{}^{1/p}}{(\rho_2-\rho_1)^2}\,\Phi(p,\rho_2). 
\end{equation}
Indeed, since by (4.23) in Lemma 4.5, 
we have 
\begin{align}
\bigg\{
\int_{B_{\rho_1}(x_0)}
\vert\tau(\varphi)\vert^{p\,\overline{\gamma}}\,v_g
\bigg\}^{1/(p\,\overline{\gamma})}
&\leq
\left(
\frac{p}{p-1}
\right)^{2/p}\frac{{C'}^{2/p}}{(\rho_2-\rho_1)^{2/p}}
\nonumber\\
&\quad\times
\bigg\{
\int_{B_{\rho_2}(x_0)}\vert 
\tau(\varphi)\vert^p\,v_g
\bigg\}^{1/p}
\end{align}
Therefore, by (5.47) and (5.58), we have 
\begin{align}
\Phi(p\,\overline{\gamma},\rho_1)
&=\bigg\{
\int_{B_{\rho_1}(x_0)}
\vert d\varphi\vert^{p\,\overline{\gamma}}\,v_g
\bigg\}^{1/{(p\,\overline{\gamma})}}
+
\bigg\{
\int_{B_{\rho_1}(x_0)}
\vert\tau(\varphi)\vert^{p\,\overline{\gamma}}\,v_g
\bigg\}^{1/{(p\,\overline{\gamma})}}
+1\nonumber\\
&\leq 
\frac{\max(p^{2/p},\,p^{m/(2p)})\,{D''}^{1/p}}
{
(\rho_2-\rho_1)^{2/p}}
\,\bigg\{
\int_{B_{\rho_2}(x_0)}\vert d\varphi\vert^p\,v_g
\bigg\}^{1/p}\nonumber\\
&\quad +
\frac{p^{1/p}\,{D''}^{1/p}}{(\rho_2-\rho_1)^{1/p}}
\bigg\{
\int_{B_{\rho_2}(x_0)}
\vert\tau(\varphi)\vert^p\,v_g
\bigg\}^{1/p}
\nonumber\\
&\quad
+
\left(
\frac{p}{p-1}
\right)^{2/p}\frac{{C'}^{2/p}}{(\rho_2-\rho_1)^{2/p}}
\bigg\{
\int_{B_{\rho_2}(x_0)}\vert 
\tau(\varphi)\vert^p\,v_g
\bigg\}^{1/p}
+1
\nonumber\\
&\leq 
\frac{\max(p^{2/p},\,p^{m/(2p)})\,D_4{}^{1/p}}{(\rho_2-\rho_1)^{2/p}}
\,\Phi(p,\rho_2), 
\end{align}
where $D_4=\max(D'', \,C')$. 
Namely, we have (5.57). 
\par
({\it The third step}))\quad 
For $k=1,2,\cdots$, let 
$p_k=m\,\overline{\gamma}^{k-1}$, 
and 
$r_k=(1+\frac{1}{2^{k-1}})\,r$. 
Then, 
$r_k$ goes to $r$ as $k\rightarrow\infty$, 
and 
apply (5.57) for 
$p=p_k$, $\rho_1=r_k<\rho_2=r_{k-1}$. 
Notice that 
$\rho_2-\rho_1=r_{k-1}-r_k=\frac{r}{2^{k-1}}$. 
Then, we have by (5.57), 
\begin{align}
\Phi(p_k,r_k)&=\Phi(p_{k-1}\,
\overline{\gamma},r_k)\nonumber\\
&\leq
\frac{
\max(p_{k-1}{}^{2/p_{k-1}},
\,p_{k-1}{}^{m/(2p_{k-1})})
\,D_4{}^{1/p_{k-1}}
}
{(\rho_2-\rho_1)^{2/p_{k-1}}}
\,\Phi(p_{k-1},r_{k-1})\nonumber\\
&
=\frac{(p_{k-1})^{m'/p_{k-1}}\,D_4{}^{1/p_{k-1}}}{(\rho_2-\rho_1)^{2/p_{k-1}}}\,\Phi(p_{k-1},r_{k-1})
\nonumber\\
&=\frac{
(m\,\overline{\gamma}^{k-2})
^{m'/(m\,\overline{\gamma}^{k-2})}
\,D_4{}^{1/(m\,\overline{\gamma}^{k-2})}
}
{
(2^{-(k-1)})^{2/(m\,\overline{\gamma}^{k-2})}
\,r^{2/(m\,\overline{\gamma}^{k-2})}
}
\,\Phi(p_{k-1},r_{k-1})\nonumber\\
&=
\frac{
C_k{}^{(k-1)/(m\,\overline{\gamma}^{k-2})}
}
{
r^{2/(m\,\overline{\gamma}^{k-2})}
}
\,\Phi(p_{k-1},r_{k-1}), 
\end{align}
where $m':=\max(2,\,\frac{m}{2})$, $\overline{\gamma}=\frac{m}{m-2}$, and 
\begin{equation}
C_k:=
(m^{m'}D_4\,\overline{\gamma}^{-m'})^{1/(k-1)}
(\overline{\gamma}^{m'}2^2)
=(D_4\,(m-2)^{m'})^{1/(k-1)}(\overline{\gamma}^{m'}2^2).
\end{equation}
Since $D_4\,(m-2)^{m'}>1$, we note that 
\begin{equation}
\overline{\gamma}^{m'}2^2
\leq\cdots\leq
C_k\leq C_{k-1}\leq\cdots\leq C_3\leq C_2
=D_4(m-2)^{m'}(\overline{\gamma}^{m'}2^2).
\end{equation}
\par
({\it The fourth step})\quad 
Therefore, we can carry out the iteration (5.60), (5.61) and (5.62). We have 
\begin{align}
\Phi(p_k,r_k)&\leq
\frac{
C_k{}^{(k-1)/(m\,\overline{\gamma}^{k-2})}
}
{
r^{2/(m\,\overline{\gamma}^{k-2})}
}
\,\Phi(p_{k-1},r_{k-1})\nonumber\\
&\leq
\frac{
C_k{}^{(k-1)/(m\,\overline{\gamma}^{k-2})}
}
{
r^{2/(m\,\overline{\gamma}^{k-2})}
}
\,
\frac{
C_{k-1}{}^{(k-2)/(m\,\overline{\gamma}^{k-3})}
}
{
r^{2/(m\,\overline{\gamma}^{k-3})}
}\,
\Phi(p_{k-2},r_{k-2})\nonumber\\
&\leq\cdots\nonumber\\
&\leq
\frac{
C_2{}^{(k-1)/(m\,\overline{\gamma}^{k-2})+
(k-2)/(m\,\overline{\gamma}^{k-3})
+\cdots+1/m
}
}
{
r^{2/(m\,\overline{\gamma}^{k-2})+
2/(m\,\overline{\gamma}^{k-3})+\cdots+
2/m
}
}
\,\Phi(p_1,r_1). 
\end{align}
Here, we have 
\begin{align}
&\frac{k-1}{m\,\overline{\gamma}^{k-2}}+
\frac{k-2}{m\,\overline{\gamma}^{k-3}}
+\cdots+\frac{1}{m}
\leq
\frac{1}{m}
\sum_{\ell=1}^{\infty}
\frac{\ell}{\overline{\gamma}^{\ell-1}}<\infty,\\
&
\frac{2}{m\,\overline{\gamma}^{k-2}}+
\frac{2}{m\,\overline{\gamma}^{k-3}}+\cdots+
\frac{2}{m}=
\frac{2}{m}\,\frac{1-(\frac{1}{\overline{\gamma}})^{k-1}}{1-\frac{1}{\overline{\gamma}}}
=
1-\big(\frac{1}{\overline{\gamma}}\big)^{k-1}
\end{align}
since 
$1<\overline{\gamma}=\frac{m}{m-2}$
and $1-\frac{1}{\overline{\gamma}}=\frac{2}{m}$. 
Thus, for the right hand side of (5.63), 
the coefficient of $\Phi(p_1,r_1)$  is estimated as 
\begin{equation}
\frac{
C_2{}^{(k-1)/(m\,\overline{\gamma}^{k-2})+
(k-2)/(m\,\overline{\gamma}^{k-3})
+\cdots+1/m
}
}
{
r^{2/(m\,\overline{\gamma}^{k-2})+
2/(m\,\overline{\gamma}^{k-3})+\cdots+
2/m
}
}
\leq
\frac{C^{\ast}}{r^{
1-(\frac{1}{\overline{\gamma}})^{k-1}
}},
\end{equation}
where
$
C^{\ast}:=C_2{}^{
\frac{1}{m}
\sum_{\ell=1}^{\infty}
\frac{\ell}{\overline{\gamma}^{\ell-1}}
}$ is a positive constant. 
\par
({\it The fifth step})\quad 
Finally, by tending $k$ to $\infty$, 
$p_k=m\,\overline{\gamma}^{k-1}\rightarrow\infty$ 
since $\overline{\gamma}>1$, and 
$r_k=(1+\frac{1}{2^{k-1}})\,r\rightarrow r$, we obtain, 
due to (5.63), 
\begin{equation}
\Phi(\infty, r)\leq \frac{C^{\ast}}{r}\,\Phi(p_1,r_1). 
\end{equation}
However, by definition (5.56) of $\Phi(p,\rho)$, 
  we have 
  \begin{equation}
  \Phi(\infty,r)=\sup_{B_r(x_0)}\vert d\varphi\vert
  +\sup_{B_r(x_0)}\vert\tau(\varphi)\vert+1,
  \end{equation}
  and 
  \begin{align}
\Phi(p_1,r_1)&=
\Phi(m,2r)\nonumber\\
&=\bigg\{
\int_{B_{2r}(x_0)}\vert d\varphi\vert^m\,v_g
\bigg\}^{1/m}
+\bigg\{
\int_{B_{2r}(x_0)}\vert \tau(\varphi)\vert^m\,v_g
\bigg\}^{1/m}+1\nonumber\\
&\leq
\varepsilon_1{}^{1/m}+\varepsilon_2{}^{1/m}+1.
  \end{align}
 Thus, we have that 
  \begin{equation}
  \sup_{B_r(x_0)}\vert d\varphi\vert
  +\sup_{B_r(x_0)}\vert\tau(\varphi)\vert+1
  \leq\frac{C^{\ast}}{r}
  \{\varepsilon_1{}^{1/m}+\varepsilon_2{}^{1/m}+1
  \}. 
  \end{equation}
  We obtain Theorem 5.1.
  \qed
\section{Proof of Theorem 3.1}
Now we are in the position to give a proof of Theorem 3.1. 
\par
Take any sequence 
$\{\varphi_i\}$ in $\mathcal F$. 
For the $\varepsilon_0>0$ in Proposition 4.3 
and for $\varepsilon_1>0$ in Theorem 5.1, 
we set 
$$
\varepsilon^{\ast}:=\min\{\varepsilon_0,\,\varepsilon_1\}>0. 
$$
Let us consider 
\begin{equation}
{\mathcal S}
:=
\left\{
x\in M\vert\,
\liminf_{i\rightarrow\infty}
\int_{B_r(x)}
\vert d\varphi_i\vert^m\,v_g\geq\varepsilon^{\ast}
\quad (\text{for all}\,\,\,r>0)
\right\}.
\end{equation}
Then, 
the set 
${\mathcal S}$ is finite. 
To see this, 
for every finite subset 
$\{x_i\}_{i=1}^k$ in $\mathcal S$, 
let us take a sufficiently small positive number 
$r_0>0$ in such a way that 
$B_{r_0}(x_i)\cap B_{r_0}(x_j)=\emptyset 
\quad (i\not=j)$. Then, we have 
for a sufficiently large $i$, 
\begin{align}
k\,\varepsilon^{\ast}&\leq
\sum_{j=1}^k\int_{B_{r_0}(x_j)}\vert d\varphi_i\vert^m\,v_g\nonumber\\
&=
\int_{\cup_{j=1}^k B_{r_0}(x_j)}\vert d\varphi_i\vert^m\,v_g\nonumber\\
&\leq
\int_M\vert d\varphi_i\vert^m\,v_g\nonumber\\
&\leq C<\infty
\end{align}
by definition of $\mathcal F$. Thus, we have 
$$
k\leq 
\frac{C}{\varepsilon^{\ast}},
$$
which implies that 
$\# {\mathcal S}\leq \frac{C}{\varepsilon^{\ast}}<\infty$. 
\par
Then, by taking a subsequence of 
$\{\varphi_i\}$ if necessary, we may assume that 
\begin{equation}
{\mathcal S}=
\left\{
x\in M\vert\,
\limsup_{i\rightarrow\infty}
\int_{B_r(x)}
\vert d\varphi_i\vert^m\,v_g\geq\varepsilon^{\ast}
\right\}.
\end{equation}
For, 
if not so, let us denote the right hand side of (6.3) by 
$\overline{\mathcal S}$. 
Then, by definition, ${\mathcal S}$ is a proper subset of 
$\overline{\mathcal S}$.  Take a point $\overline{x}\in \overline{\mathcal S}\backslash {\mathcal S}$. 
By taking a subsequence 
of $\{\varphi_i\}$, 
by the same letter, in such a way that 
$$
\liminf_{i\rightarrow \infty}\int_{B_r(\overline{x})}
\vert d\varphi_i\vert^m\,v_g
\geq\varepsilon^{\ast}, 
$$
For this $\{\varphi_i\}$, $\overline{x}$ belongs to 
$\mathcal S$. Since $\mathcal S$ is a finite set, 
this process stops at finite times, then at last we have 
$\overline{\mathcal S}={\mathcal S}$. 
\par
Now, let 
$x\in M\backslash{\mathcal S}$. 
Then, 
\begin{equation}
\limsup_{i\rightarrow\infty}\int_{B_r(x)}\vert d\varphi_i\vert^m\,v_g<\varepsilon^{\ast}. 
\end{equation}
Due to Proposition 4.3 and the definition of $\mathcal F$, 
we have
\begin{align}
\sup_{B_{r}(x)}\vert\tau(\varphi_i)\vert^2
&\leq
\frac{C}{r^{m/2}}\int_{B_{2r}(x)}\vert\tau(\varphi_i)\vert^2\,v_g
\nonumber\\
&\leq 
\frac{C^2}{r^{m/2}}, 
\end{align}
so that we have that 
\vskip0.3cm\par
\qquad 
{\it $(C^0):$ \qquad
 the $C^0$-estimate on $B_r(x)$ 
of $\tau(\varphi_i)$ 
uniformly on $i$. }
\vskip0.3cm\par
And due to (6.5), for a sufficiently small 
$r>0$, we can show that 
\begin{equation}
\int_{B_{r/2}(x)}\vert\tau(\varphi_i)\vert^m\,v_g
<\varepsilon_2,
\end{equation}
where $\varepsilon_2>0$ is the constant in 
Theorem 5.1. 
Indeed, we have 
due to (4.7) in Proposition 4.3, 
\begin{align}
\int_{B_{r/2}(x)}\vert\tau(\varphi_i)\vert^m\,v_g
&\leq
\sup_{B_{r/2}(x)}\vert\tau(\varphi_i)\vert^{m-2}
\int_{B_{r/2}(x)}\vert\tau(\varphi_i)\vert^2\,v_g
\nonumber\\
&=
(\sup_{B_{r/2}(x)}\vert\tau(\varphi_i)\vert^2
)^{\frac{m-2}{2}}\,
\int_{B_{r/2}(x)}\vert\tau(\varphi_i)\vert^2\,v_g
\nonumber\\
&\leq
(\sup_{B_{r/2}(x)}\vert\tau(\varphi_i)\vert^2
)^{\frac{m-2}{2}}\,
\int_{B_{r}(x)}\vert\tau(\varphi_i)\vert^2\,v_g
\nonumber\\
&\leq
\bigg(
\frac{C'}{r^{m/2}}\int_{B_r(x)}\vert\tau(\varphi_i)\vert^2\,v_g
\bigg)^{\frac{m-2}{2}}\,\int_{B_r(x)}\vert\tau(\varphi_i)\vert^2\,v_g
\nonumber\\
&=\frac{(C')^{\frac{m-2}{2}}}{r^{\frac{m}{2}\,\frac{m-2}{2}}}\,\bigg(
\int_{B_r(x)}\vert\tau(\varphi_i)\vert^2\,v_g
\bigg)^{\frac{m}{2}}.
\end{align}
Here,  
substituting the inequality (6.5) into 
$\int_{B_r(x)}\vert\tau(\varphi_i)\vert^2\,v_g$, 
we have 
\begin{align}
\int_{B_r(x)}\vert\tau(\varphi_i)\vert^2\,v_g
&\leq 
\sup_{B_r(x)}\vert\tau(\varphi_i)\vert^2\,\text{Vol}(B_r(x))
\nonumber\\
&\leq\frac{C^2}{r^{\frac{m}{2}}}\,\,
C_0r^m
=C_0C^2\,\,r^{\frac{m}{2}}
\end{align}
since Vol$(B_r(x))\leq \alpha_{m-1}\int^r_0
\frac{1}{\sqrt{-\delta}}\sinh^{m-1}(\sqrt{-\delta}\,t)\,dt\leq C_0\,r^m$ for a sufficiently small $r>0$ for some constant $C_0>0$, where $\alpha_{m-1}$ 
is the volume of the $(m-1)$-dimensional unit sphere. 
This follows from that the Ricci curvature 
of $(M,g)$ satisfies 
$\text{\rm Ric}^M\geq (m-1)\delta$ 
for some negative number $\delta$ since $M$ is compact. 
Therefore, by (6.8), 
the last line of the right hand side of (6.7) is estimated as 
\begin{align}
\frac{(C')^{\frac{m-2}{2}}}{r^{\frac{m}{2}\,\frac{m-2}{2}}}\,\bigg(
\int_{B_r(x)}\vert\tau(\varphi_i)\vert^2\,v_g
\bigg)^{\frac{m}{2}}&\leq 
\frac{(C')^{\frac{m-2}{2}}}{r^{\frac{m}{2}\,\frac{m-2}{2}}}\,
\big(
C_0C^2\,r^{\frac{m}{2}}
\big)^{\frac{m}{2}}\nonumber\\
&=(C')^{\frac{m-2}{2}}\,(C_0C^2)^{\frac{m}{2}}\,
r^{\frac{m}{2}}
\end{align}
which tends to zero if $r\rightarrow 0$. So, we have (6.6). Notice that we may regard 
this $\frac{r}{2}$, to be $r$ in (6.6). 
\par
Due to the inequality $\varepsilon^{\ast}\leq \varepsilon_1$ where $\varepsilon_1$ is the constant in Theorem 5.1, 
we can apply Theorem 5.1 to $\varphi_i$, which implies that 
\begin{equation}
\sup_{B_{r/2}(x)}\vert d\varphi_i\vert+
\sup_{B_{r/2}(x)}\vert\tau(\varphi_i)\vert
\leq 
\frac{C^{\ast}}{r}\,(\varepsilon_1{}^{1/m}+
\varepsilon_2{}^{1/m}+1),
\end{equation} 
so we have that 
\vskip0.3cm\par\quad
{\it $(C^1)$\qquad the $C^1$-estimate 
on $B_r(x)$ of $\varphi_i$ uniformly on $i$.} 
\vskip0.3cm\par
On the other hand, since
$\varphi_i\in {\mathcal F}$, 
all the $\varphi_i$ are biharmonic, i.e., 
they satisfy the equations 
\begin{align}
&\tau_2(\varphi_i)=\overline{\Delta}(\tau(\varphi_i))
-{\mathcal R}(\tau(\varphi_i))=0\\
&\Longleftrightarrow
\left\{
\begin{aligned}
&(1)\quad \overline{\Delta}\sigma_i={\mathcal R}(\sigma_i),\\
&(2)\quad \tau(\varphi_i)=\sigma_i.
\end{aligned}
\right.
\end{align}
\par
We have
\begin{equation*}
{\mathcal R}(\sigma_i)=
\sum_{j=1}^m{}^N\!R(\sigma_i,d\varphi_i(e_j))d\varphi_i(e_j), 
\end{equation*}
and $\vert d\varphi_i(e_j)\vert$ $(j=1,\cdots,m)$ is bounded 
uniformly on $i$ since the $C^1$-estimate on 
$B_r(x)$ of $\varphi_i$ uniformly on $i$. 
Due to (1) of (6.12), 
each $\sigma_i$ are solutions of 
the linear equations with coefficients 
which are bounded uniformly on $i$. 
Due to the standard argument 
(cf. Theorem 3.1, p.397 in  \cite{LU}), 
we have 
\vskip0.3cm\par
 \quad {\it $(C^{\alpha})$\qquad 
 the $C^{\alpha}$-estimate on $B_r(x)$ 
 of $\sigma_i$ uniformly on $i$. }
\vskip0.3cm \par
 Furthermore, the equation (2) of (6.12) means that $\varphi_i$ are solutions of the non-linear equations 
 \begin{align}
 \tau(\varphi_i)^{\gamma}=
 \Delta(\varphi_i{}^{\gamma})+
 \sum_{j, k=1}^m\sum_{\alpha,\beta=1}^n
 g^{jk}\,
 {}^N\!\Gamma^{\gamma}_{\alpha\beta}(\varphi_i)\frac{\partial\varphi_i{}^{\alpha}}{\partial x_j}
 \frac{\partial\varphi_i{}^{\beta}}{\partial x_k}
 =\sigma_i{}^{\gamma},
 \end{align}
where $(x_1,\cdots,x_m)$, and $(y_1,\cdots,y_n)$ are the coordinates of $M$, and $N$, respectively, and we denote $\varphi_i{}^{\alpha}
 =y_{\alpha}\circ\varphi_i$. 
 Namely, we simply denote (6.13) by 
 \begin{align}
 \tau(\varphi_i)=\Delta\varphi_i
 +{}^N\!\Gamma(\varphi_i)(d\varphi_i,d\varphi_i)=\sigma_i,
 \end{align}
 with bounded coefficients on $\varphi_i$ since
 $\vert d\varphi_i\vert$ are bounded due to the $C^1$-estimate of $\varphi_i$. 
 Notice here that we regard in the second term of (6.14) one of two $d\varphi_i$ should be coefficients, and then $\varphi_i$ can be regarded to be solutions of the linear equations. 
 Then, due to the standard argument (cf. Theorem 4.1, p. 399 in \cite{LU}) 
 for (6.14), we have 
 \vskip0.3cm\par\quad
 $(C^{1,\alpha})$ \qquad 
 {\it the $C^{1,\alpha}$-estimate on $B_r(x)$ 
 of $\varphi_i$ uniformly on $i$.
 }
 \vskip0.3cm\par
 Hence, we can regard again to be solutions of the linear equation (6.14) with the $C^{\alpha}$-coefficient $d\varphi_i$. Then, due to the Schauder estimate, we have
 \begin{align}
 \vert \varphi_i\vert_{C^{2,\alpha}(B_r(x))}
 &\leq
 C(\vert\varphi_i\vert_{C^0(B_r(x))}+\vert\sigma_i\vert_{C^{\alpha}(B_r(x))})
 \nonumber\\
 &\leq C(C_1+C_2), 
 \end{align}
 because that  
 $\vert\varphi_i\vert\leq C_1$ follows from the $C^0$-estimate of $\varphi_i$ and that 
 $\vert\sigma_i\vert\leq C_2$ follows from 
 the $C^{\alpha}$-estimate of $\sigma_i$. 
 \par Thus, we have 
\vskip0.3cm\par \quad 
 $(C^{2,\alpha})$\qquad the $C^{2,\alpha}$-estimate on $B_r(x)$ of 
 $\varphi_i$ uniformly on $i$. 
 \vskip0.3cm\par
Finally, due to the bootstrap argument, 
we obtain the $C^{\infty}$-estimate on 
$B_r(x)$ 
of $\varphi_i$ uniformly on $i$.  
\vskip0.3cm\par
Therefore, there exists 
a subsequence $\{\varphi_{i_j}\}$ of 
$\{\varphi_i\}$ 
and a smooth map 
$\varphi_{\infty}:\,M\backslash{\mathcal S}\rightarrow N$ 
such that 
$\varphi_{i_j}$ converges to $\varphi_{\infty}$ on 
$B_r(x)$ in the $C^{\infty}$-topology 
as $j\rightarrow\infty$. Thus, $\varphi_{\infty}:(M\backslash{\mathcal S},g)\rightarrow(N,h)$ is also biharmonic. 
\par
For (2) in Theorem 3.1, 
let us consider the Radon measures 
$\vert d\varphi_{i_j}\vert^m\,v_g$. 
Then, these have a weak limit which is also a 
Randon measure, say $\mu$. 
Recall that $\mu$ is by definition a {\em Radon measure} if 
$(1)$ $\mu$ is locally finite, i.e., 
$\mu(K)<\infty$ for every 
compact subset $K$ on $M$, 
and 
$(2)$ $\mu$ is Borel regular, i.e., 
it holds that, for all Borel subset $A$ of $M$, 
\begin{align*}
\mu(A)&=\sup\{\mu(K)\vert\,\,
\text{for all compact subset $K$ of $A$}\},\,\,\text{and}\\
\mu(A)&=\inf\{\mu(O)\vert\,\,
\text{for all open subset $O$ of $M$  including $A$}\}. 
\end{align*}
On the other hand, since $\varphi_{i_j}$ converges 
to 
$\varphi_{\infty}$ 
on $M\backslash{\mathcal S}$ in the $C^{\infty}$-topology 
as $j\rightarrow\infty$, 
it holds that 
\begin{equation}
\mu=\vert d\varphi_{\infty}\vert^m\,v_g \quad \text{on}\,\,
M\backslash {\mathcal S}.
\end{equation} 
Here, $\mathcal S$ is a finite subset of $M$, say 
${\mathcal S}=\{x_1,\cdots,x_k\}$. 
Then, the Radon measure $\mu-\vert d\varphi_{\infty}\vert^m\,v_g$ satisfies that its support is contained in $\mathcal S$. 
Therefore, 
it holds that 
\begin{equation}
\mu-\vert d\varphi_{\infty}\vert^m\,v_g=
\sum_{i=1}^ka_j\,\delta_{x_j}
\end{equation}
for some non-negative real numbers $a_j$ $(j=1,\cdots,k)$.  
Remark here that 
$a_j<\infty$ for every $j=1,\cdots, k$. 
Because $\mu$ is a Radon measure, so that $\mu$ is locally finite. 
Therefore, 
the Radon measure $\vert d \varphi_{i_j}\vert^m\,v_g$ 
converges weakly to 
$\mu$, and 
\begin{equation}
\mu=\vert d \varphi_{\infty}\vert^m\,v_g+\sum_{j=1}^ka_j\,\delta_{x_j}
\end{equation}
due to (6.9). 
We have (2) of Theorem 3.1. 
\qed
\vskip1.4cm\par

\end{document}